\documentclass[11pt]{amsart}
\usepackage{amscd,amssymb,subfigure}
\usepackage[arrow,matrix,graph]{xy}
\usepackage{graphicx}
\usepackage[colorlinks,plainpages,backref,urlcolor=blue]{hyperref}

\newtheorem{theorem}{Theorem}[section]
\newtheorem{corollary}[theorem]{Corollary}
\newtheorem{lemma}[theorem]{Lemma}
\newtheorem{prop}[theorem]{Proposition}

\newtheorem{thm}{Theorem}

\theoremstyle{definition}
\newtheorem{definition}[theorem]{Definition}
\newtheorem{example}[theorem]{Example}
\newtheorem{remark}[theorem]{Remark}
\newtheorem*{ack}{Acknowledgment}

\newcommand{\N}{\mathbb{N}}
\newcommand{\Z}{\mathbb{Z}}
\newcommand{\Q}{\mathbb{Q}}

\newcommand{\C}{\mathbb{C}}

\newcommand{\T}{\mathbb{T}}
\newcommand{\K}{\mathbb{K}}
\newcommand{\PP}{\mathbb{P}}

\renewcommand{\AA}{\mathsf{A}}

\newcommand{\p}{{\mathfrak p}}
\newcommand{\m}{{\mathfrak m}}

\newcommand{\V}{{\mathcal V}}
\newcommand{\E}{{\mathcal E}}
\newcommand{\LL}{{\mathcal L}}
\newcommand{\OO}{{\mathcal O}}

\newcommand{\A}{{\mathcal A}}

\DeclareMathOperator{\rank}{rank}

\DeclareMathOperator{\im}{im}
\DeclareMathOperator{\Hom}{Hom}
\DeclareMathOperator{\coker}{coker}
\DeclareMathOperator{\codim}{codim}
\DeclareMathOperator{\id}{id}
\DeclareMathOperator{\ab}{{ab}}
\DeclareMathOperator{\abf}{{abf}}

\DeclareMathOperator{\Sp}{Sp}

\DeclareMathOperator{\ann}{{ann}}
\DeclareMathOperator{\Spec}{{Spec}}

\DeclareMathOperator{\lk}{{lk}}

\DeclareMathOperator{\df}{def}

\DeclareMathOperator{\dire}{dir}
\DeclareMathOperator{\Tors}{Torsion}

\newcommand{\surj}{\twoheadrightarrow}
\newcommand{\inj}{\hookrightarrow}
\newcommand{\isom}{\xrightarrow{\,\simeq\,}}
\newcommand{\twoarr}{\Longleftrightarrow}
\newcommand{\abs}[1]{\left| #1 \right|}

\newcommand{\vv}{{\mathcal W}_1}
\newcommand{\w}{W_1}
\newcommand{\bg}{b_1^{\rm gen}}
\newcommand{\const}{\rm const}

\newenvironment{romenum}
{

\begin{enumerate}}{\end{enumerate}}

\begin{document}
%\date{June 16, 2007}

\title[Alexander polynomials: Essential variables and multiplicities]%
{Alexander polynomials: \\ Essential variables and multiplicities}

\author[A.~Dimca]{Alexandru Dimca}
\address{Laboratoire J.A.~Dieudonn\'{e}, UMR du CNRS 6621, 
Universit\'{e} de Nice--Sophia Antipolis, Parc Valrose,
06108 Nice Cedex 02, France}
\email{dimca@math.unice.fr}

\author[S.~Papadima]{Stefan Papadima$^1$}
\address{Inst.~of Math.~Simion Stoilow, 
P.O. Box 1-764,
RO-014700 Bucharest, Romania}
\email{Stefan.Papadima@imar.ro}

\author[A.~Suciu]{Alexander~I.~Suciu$^2$}
\address{Department of Mathematics,
Northeastern University,
Boston, MA 02115, USA}
\email{a.suciu@neu.edu}
%\urladdr{http://www.math.neu.edu/\~{}suciu}

\thanks{$^1$Partially supported by the CEEX Programme of
the Romanian Ministry of Education and Research, contract
2-CEx 06-11-20/2006.}
\thanks{$^2$Partially supported by NSF grant DMS-0311142}



\begin{abstract}
We explore the codimension one strata in the degree-one 
cohomology jumping loci of a finitely generated group, through   
the prism of the multivariable Alexander polynomial. As an 
application, we give new criteria that must be satisfied by 
fundamental groups of smooth, quasi-projective complex varieties. 
These criteria establish precisely which fundamental groups of 
boundary manifolds of complex line arrangements 
are quasi-projective. We also give sharp upper bounds for 
the twisted Betti ranks of a group, in terms of  multiplicities 
constructed from the Alexander polynomial. For Seifert links in 
homology $3$-spheres, these bounds become equalities, and 
our formula shows explicitly how the Alexander polynomial 
determines all the characteristic varieties.
\end{abstract}
\maketitle

\tableofcontents

\section{Introduction and statement of results}
\label{sect:intro}

\subsection{Characteristic varieties}
\label{intro:1} 
Let $G$ be a finitely generated group.
Each character $\rho\colon G\to \C^*$ in the complex 
algebraic group $\T_G= \Hom (G, \C^*)$ gives rise to a 
rank~$1$ local system on $G$, denoted by $\C_{\rho}$. 
The {\em characteristic varieties} of $G$ are the 
jumping loci for homology with coefficients 
in such local systems:
\[
\V_k (G)= \{ \rho\in \T_G  \mid \dim_{\C} H_1(G; \C_{\rho})\ge k \}.
\]

An alternate description is as follows.  Let $\Z{G}$ be 
the group ring of $G$, with augmentation 
ideal $I_G$, and let $G_{\abf}$ be the maximal torsion-free 
abelian quotient of $G$. Note that $\C{G}_{\abf}$ is the 
coordinate ring of identity component of the character torus, 
$\T^0_G$. Finally, let $A_G=\Z{G_{\abf}} \otimes_{\Z{G}} I_G$ 
be the Alexander module. Then, at least away from the origin, 
$\V_k (G) \cap \T^0_G$ coincides with the subvariety of 
$\T^0_G$ defined by $\E_k (A_G)$, the ideal of codimension 
$k$ minors of a presentation matrix for $A_G$. 

We focus in this paper on $\vv (G)$, the union of all codimension-one 
irreducible components of $\V_1(G)$, which are contained in 
$\T^0_G$. It turns out that this variety is closely related 
to another, even more classical invariant. 

\subsection{The Alexander polynomial}
\label{intro:2}
The group ring $\Z{G}_{\abf}$ may be identified with the 
Laurent polynomial ring 
$\Lambda_n=\Z[t_1^{\pm 1},\dots ,t_n^{\pm 1}]$, 
where $n=b_1(G)$.  Our other object of study is the 
(multivariable) {\em Alexander polynomial},  $\Delta^G$, 
defined as the greatest common divisor of all elements of 
$\E_1(A_G)$, up to units in $\Lambda_n$. 

After reviewing these basic notions in \S\ref{sect:prel}, in a
more general context, we start in \S\ref{sect:alexvc1} 
by making precise the relationship between the Alexander polynomial 
of $G$ and the codimension-one stratum of $\V_1 (G)$. For example, 
if $b_1(G)> 1$ and $\Delta^G\ne 0$, then $\vv(G)$ is simply the 
hypersurface defined by $\Delta^G$; for the remaining cases, 
see Corollary \ref{cor:deltacd1}.  As an application, we compute  
in Proposition \ref{prop:delta arr} the Alexander polynomials of 
arrangement groups: if $\A$ is a finite set of hyperplanes in some 
$\C^{\ell}$, and $G$ is the fundamental group of its complement, 
then $\Delta^G$ is constant, except when $\A$ is a pencil, 
and $\abs{\A}\ge 3$.  

We then ask: When does the Alexander polynomial 
of $G$ have a {\em single essential variable}, i.e., 
\[
\Delta^G (t_1,\dots , t_n)\doteq  P(t_1^{e_1}\cdots t_n^{e_n}),
\] 
for some polynomial $P\in \Z[t^{\pm 1}]$?  Equivalently, when 
is the Newton polytope of $\Delta^G$ a line segment? 
If $b_1(G)\ge 2$, this happens precisely when $\vv (G)=\emptyset$, 
or all the irreducible components of $\vv (G)$ are parallel, 
codimension-one subtori of $\T_G^0$; see 
Proposition \ref{prop:cod1ess}.

\subsection{Obstructions to quasi-projectivity}
\label{intro:3}
Some 50 years ago, work of J.-P. Serre \cite{Se} raised 
the following problem: characterize those finitely presented 
groups that appear as the fundamental groups of smooth, 
connected, (quasi-) projective complex varieties. 
Such groups are called {\em (quasi-) projective groups}, 
see e.g.~\cite{Ar95}, \cite{Ca03}.

A major goal in this note is to present and exploit new 
obstructions to quasi-projec\-tivity of groups.  Our 
approach in \S\ref{sect:qprojalex} is based on a 
foundational result of Arapura \cite{A}, which states: 
If $W$ is an irreducible component of the first characteristic 
variety of a quasi-projective group $G$, then 
$W=\rho \cdot \dire W $, with $\rho\in \T_G$ 
and $\dire W\subset \T_G^0$  a subtorus through $1$.  
Using refinements from \cite{DPS} and \cite{D06}, we formulate 
our first set of quasi-projectivity obstructions in terms of the 
relative position of the components of $\V_1(G)$, as follows 
(a more detailed statement is given in  
Theorem \ref{thm:posobs}). 

\begin{thm}
\label{thmA}
Let $G$ be a quasi-projective group. 
If $W$ and $W'$ are two distinct components of $\V_1(G)$, 
then either $\dire W= \dire W'$, or 
$\dire W \cap \dire W'$ is a finite set.
\end{thm}

Using the dictionary between $\vv(G)$ and $\Delta^G$, we find that, 
with one exception, the Alexander polynomial of a 
quasi-projective group must have a single essential variable. 
More precisely, we prove the following result (a more detailed 
statement is given in Theorem \ref{thm:alexone}). 

\begin{thm}
\label{thmB}
Let $G=\pi_1(M)$ be the fundamental group of a smooth, connected, 
complex quasi-projective variety. If $b_1(G)\ne 2$, then 
the Alexander polynomial $\Delta^G$ has a single essential variable. 
If $M$ is projective, then $\Delta^G \doteq \const$. 
\end{thm}

The condition on the first Betti number of $G$ is really necessary. 
In Example \ref{ex:bdoi}, we exhibit a smooth, non-compact algebraic 
surface $M$, with fundamental group $G$, for which $b_1(G)=2$ 
and the Newton polytope of $\Delta^G$ is $2$-dimensional. 

As explained in Example \ref{ex:alexvsgen}, the characteristic varieties
and Alexander polynomials of finitely presented groups can be very complicated.
The restrictions imposed by our theorems are quite efficient 
when attacking Serre's problem. As an illustration, we consider 
groups of the form $G_{\A}= \pi_1(M_{\A})$, where $\A$ is an 
arrangement of lines in the complex projective plane $\PP^2$, 
and $M_{\A}$ is the boundary of a regular neighborhood 
of $\A$. Using an explicit formula for 
$\Delta^{G_{\A}}$, recently obtained in \cite{CS06}, we 
determine in Proposition \ref{prop:serrebd} the precise class 
of arrangements $\A$ for which $G_{\A}$ is a quasi-projective group.

\subsection{Twisted Betti ranks}
\label{intro:4}
In Sections \ref{sect:multupper} and \ref{sect:aai}, 
we explore the connection between the multiplicities 
of the factors of the Alexander polynomial, $\Delta^G$, 
and the higher-depth characteristic varieties, $\V_k(G)$.
We use this connection to give easily computable, 
sharp upper bounds for the `twisted' Betti ranks,
$b_1(G, \rho)= \dim_{\C} H_1(G; \C_{\rho})$, 
corresponding to non-trivial characters $\rho\in \T_G$. 

Before stating our results, we need some more terminology. 
Given an irreducible Laurent polynomial $f\in \C{G}_{\abf}$, define 
the `generic' Betti number, $\bg (G, f)$, as  the largest integer 
$k$ so that $\V_k(G)$ contains $V(f)$. It turns out that 
$b_1(G, \rho)=\bg (G,f)$, for $\rho$ in a certain non-empty 
Zariski open subset of $V(f)$. Also, denote by $\nu_{\rho}(f)$ the order 
of vanishing of the germ of $f$ at $\rho$.  

Finally, declare the Alexander ideal $\E_1(A_G)$ to be 
{\em almost principal} if there is an integer $d\ge 0$ such that 
$I_{G_{\abf}}^d\cdot (\Delta^G)\subset \E_1(A_G)$, 
over $\C$. Examples abound: groups with $b_1(G)=1$, link groups, 
as well as fundamental groups of closed, orientable 
$3$-manifolds, see \cite{MM}. 

\begin{thm}
\label{thmC}
Let $\Delta^G \doteq_{\C} \prod_{j=1}^s f_j^{\mu_j}$ be the 
decomposition into irreducible factors of the Alexander polynomial
of a finitely generated group $G$. Then:
\begin{romenum}
\item \label{i}
$\bg (G, f_j)\le \mu_j$, for each $j=1,\dots, s$.
\item \label{ii}
If the Alexander ideal of $G$ is almost principal, 
then, for all $\rho\in \T_G^0\setminus \{1\}$, 
\[
b_1(G, \rho)\le \sum_{j=1}^s \mu_j\cdot \nu_{\rho}(f_j).
\]
\end{romenum}
\end{thm}

This result is proved in Theorems \ref{thm:uppergen} 
and \ref{thm:uppergral}. We also give a variety of examples 
showing that the inequalities can be strict, yet in general 
are sharp. When $G_{\ab}$  has no torsion and the upper bounds 
from Part \eqref{ii} above are attained, for every $\rho$, it
follows that $\Delta^G$ determines $\V_k (G)$, for all $k$.

In the setup from Theorem \ref{thmC}\eqref{ii}, suppose 
also that $\Delta^G$ has a single essential variable. Then 
the `generic' upper bounds from \eqref{i} take the form
\begin{equation}
\label{bound}  
b_1(G, \rho)\le \mu_j, \ \text{for all $\rho\in V(f_j)\setminus \{ 1\}$},
\end{equation}
see Corollary \ref{cor:upperqp}. The simplest situation is when 
$b_1(G)=1$, in which case the upper bounds \eqref{bound} are 
all attained if and only if the monodromy of the complexified 
Alexander invariant is semisimple, see Proposition \ref{prop:princtest}. 
For example, if $G=\pi_1(M)$, where $M$ is a compact, connected 
manifold fibering over the circle with connected fibers, then 
$b_1(G)=1$ and the inequalities \eqref{bound} are all equalities 
precisely when the algebraic monodromy is semisimple and with 
no eigenvalue equal to $1$, see Corollary \ref{cor:monotest}.  

\subsection{Seifert links}
\label{subsec:intro5}
We conclude in Section \ref{sect:seifert} with a study of the 
characteristic varieties of a well known class of links in  
homology $3$-spheres: those for which the link exterior
admits a Seifert fibration.  Since any such link, $(\Sigma^3,L)$,  
has an analytic representative given by an isolated, weighted 
homogeneous singularity (see  \cite{EN}, \cite{Ne}), the corresponding 
link group, $G_L=\pi_1(\Sigma^3\setminus L)$, is quasi-projective. 
In \cite{EN}, Eisenbud and Neumann study in detail the 
Alexander-type invariants of Seifert link groups, obtaining 
a formula for the Alexander polynomial 
of $G_L$, in terms of a `splice diagram' representing $L$. 

The  Eisenbud-Neumann calculus is used in 
Theorem~\ref{thm:seifert cv} to show that the upper 
bounds from \eqref{bound} all become equalities for 
Seifert links $L$.  Consequently, the Alexander polynomial of 
$L$ completely determines the characteristic 
varieties $\V_k(G_L)$, for all $k\ge 1$.

\section{Alexander invariants and characteristic varieties}
\label{sect:prel}

We start by collecting some known facts about Alexander-type invariants 
and characteristic varieties, in a slightly more general context than the one 
outlined in the Introduction. 

\subsection{Alexander polynomials}
\label{subsec:fitt}

Let $R$ be a commutative ring with unit. Assume $R$ is 
Noetherian and a unique factorization domain.  If two elements 
$\Delta$, $\Delta'$ in $R$ generate the same principal ideal (that is,  
$\Delta =u\Delta'$, for some unit $u\in R^*$), we 
write $\Delta \doteq \Delta'$. 

Let $A$ be a finitely-generated $R$-module. Then $A$ admits 
a finite presentation, $R^m\xrightarrow{M}  R^n \to A \to 0$. 
The {\em $i$-th elementary ideal} of $A$, denoted $\E_i(A)$, 
is the ideal in $R$ generated by the minors of size $n-i$ of 
the $n \times m$ matrix $M$, with the convention that 
$\E_i(A)=R$ if  $i \ge n$, and $\E_i(A)=0$ if $n-i>m$. 
Clearly, $\E_i(A) \subset \E_{i+1}(A)$, for all $i\ge 0$.  
Of particular interest is the ideal of maximal minors, 
$\E_0(A)$, also known as the {\em order ideal}.

The $i$-th {\em Alexander polynomial} of $A$, 
denoted $\Delta_i(A)$, is a generator of the smallest principal 
ideal in $R$ containing $\E_i(A)$, that is, the greatest common
divisor of all elements of $\E_i(A)$. As such, $\Delta_i(A)$ is 
well-defined only up to units in $R$. 
Note that $\Delta_{i+1}(A)$ divides $\Delta_i(A)$, 
for all $i\ge 0$. 

\begin{example} 
\label{ex0}
Suppose 
$A=R^s \oplus R/(\lambda_1) \oplus \cdots \oplus R/(\lambda_r)$, 
where $\lambda_1, \dots, \lambda_r$ are non-zero
elements in $R$ such that $\lambda_r\!\mid\cdots \mid\!  \lambda_1$. 
Then 
\begin{equation}
\label{eq:delta pid}
\Delta_i(A)\doteq\begin{cases}
0 & \text{if $0 \le i \le s-1$},
\\
\lambda_{i-s+1} \cdots \lambda_r 
& \text{if $s\le i \le s+r-1$},
\\
1& \text{if $s+r \le i$}.
\end{cases}
\end{equation}
This computation applies to any finitely generated module
over a principal ideal domain $R$; for instance, $R=\K[t]$ 
with $\K$ a field, or $R$ a discrete valuation ring.
\end{example} 

We will be mostly interested in the case when 
$R$ is the ring of Laurent polynomials over the integers, 
$\Lambda_q=\Z[t^{\pm 1}_1, \dots ,t_q^{\pm 1}]$, 
or over the complex numbers, 
$\Lambda_q \otimes \C=\C[t^{\pm 1}_1, \dots ,t_q^{\pm 1}]$. 
We will denote by $\doteq$ equality up to units in 
$\Lambda_q$ (of the form 
$u=\pm  t_1^{\nu_1} \cdots t_q^{\nu_q}$), and by 
$\doteq_\C$ equality up to units in 
$\Lambda_q\otimes \C$ (of the form 
$u=z t_1^{\nu_1} \cdots t_q^{\nu_q}$, with $z\in \C^*$). 
In particular, we will write $\Delta \doteq \const$ if 
$\Delta\doteq c$, for some $c\in \Z$ 
(equivalently, $\Delta \doteq_{\C} \text{$0$ or $1$}$). 

Suppose $A$ is a finitely-generated module over $\Lambda_q$. 
We then have Alexander polynomials, 
$\Delta_{i}(A)\in \Lambda_q$ and
$\Delta_{i}(A\otimes \C)\in \Lambda_q \otimes \C$, 
well defined up to units in the respective rings. 
Using the unique factorization property, it is readily 
seen that $\Delta_{i}(A\otimes \C) \doteq_\C \Delta_{i}(A)$. 

\subsection{Support varieties}
\label{subsec:char var}
As before, let $A$ be a finitely-generated $R$-module. 
The {\em support} of $A$ is the reduced subscheme of 
$\Spec(R)$ defined by the order ideal, $\E_0(A)$;  we 
will denote it by $V_1(A)$.  Since
$\sqrt { \E_0(A)}= \sqrt{\ann A }$, 
this is the usual notion of support in algebraic geometry, 
based on the annihilator ideal of the module $A$. 
In particular, a prime ideal $\p \subset R$ belongs to 
$V_1(A)$ if and only if the localized module $A_\p$ is
non-zero.

More generally, the {\em $k$-th support variety} 
of $A$, denoted  $V_k(A)$, is the reduced subscheme 
of $\Spec(R)$ defined by the ideal $\E_{k-1}(A)$.  
Clearly, $\{V_k(A)\}_{k\ge 1}$ is a decreasing sequence 
of subvarieties of $\Spec(R)$; these varieties 
are invariants of the $R$-isomorphism type of $A$.  

If $\codim V_k(A)>1$, then $\Delta_{k-1}(A)\doteq 1$, so 
$\Delta_{k-1}(A)$ carries no useful information.  One 
of our guiding principles in this note is to elaborate on 
this, by explaining  how much can be extracted 
from the Alexander polynomials, in the case when 
they are non-constant.

\subsection{Alexander modules}
\label{ss:tprel}
Let $X$ be a connected CW-complex. Without loss of generality, 
we may assume that $X$ has a unique $0$-cell, call it $x_0$. 
We will make the further assumption that $X$ has finitely many 
$1$-cells, in which case the fundamental group, $G=\pi_1(X,x_0)$, 
is finitely generated.   

Let $\phi\colon G \surj H$ be a homomorphism onto an abelian group 
$H$, and let $p\colon X^{\phi} \to X$ be the corresponding Galois cover. 
Denote by $F=p^{-1}(x_0)$ the fiber of $p$ over the basepoint. 
The exact sequence of the pair $(X^{\phi}, F)$ then yields 
an exact sequence of $\Z{H}$-modules,
\begin{equation}
\label{eq:crow1}
\xymatrix{0 \ar[r]& H_1(X^{\phi};\Z) \ar[r]& H_1(X^{\phi},F;\Z) \ar[r]& 
H_0(F;\Z) \ar[r]& H_0(X^{\phi};\Z) \ar[r]& 0}.
\end{equation}

The $\Z{H}$-modules $B_{\phi}=H_1(X^{\phi};\Z)$ and 
$A_{\phi} = H_1(X^{\phi},F;\Z)$ are called the {\em Alexander invariant} 
(respectively,  the {\em Alexander module}) of $X$, relative to $\phi$. 
Identifying the kernel of $H_0(F;\Z) \to H_0(X^{\phi};\Z)$ with the 
augmentation ideal, $I_H=\ker \big(\Z H \xrightarrow{\epsilon} \Z\big)$,  
yields the  `Crowell exact sequence' of $X$, 
\begin{equation}
\label{eq:crow2}
\xymatrix{0 \ar[r]& B_{\phi} \ar[r]& A_{\phi} \ar[r]& I_H \ar[r]& 0}.
\end{equation}

By considering a classifying map $X\to K(G,1)$, it is readily 
seen that both $\Z{H}$-modules, $A_{\phi}$ and 
$B_{\phi}$, depend only on the homomorphism $\phi$.  
For example, if $\phi_{\ab}\colon G \to G_{\ab}$ is the 
abelianization map (corresponding to the maximal 
abelian cover $X^{\ab}\to X$), then 
$A_{\ab}=\Z{G_{\ab}}\otimes_{\Z{G}} I_G$ 
and $B_{\ab}=G'/G''$, with $\Z{G}_{\ab}$-module 
structure on $B_{\ab}$ determined by the extension 
\begin{equation}
\label{eq:bab}
\xymatrix{0 \ar[r]& G'/G''  \ar[r]& G/G''  \ar[r]& G/G'  \ar[r]& 0}.
\end{equation}

\subsection{Torsion-free abelian covers}
\label{subsec:stuff}
Assume now $H$ is a finitely-generated, free abelian group; 
let $q$ be its rank. Identify the group ring $\Z{H}$ with 
$\Lambda_q=\Z[t_1^{\pm 1}, \dots , t_q^{\pm 1}]$, the 
augmentation ideal $I_H$ with the ideal $I=(t_1-1,\dots , t_q-1)$,   
and $\T_H=\Spec(\C{H})$ with the algebraic 
torus $\Hom(H,\C^*)=(\C^*)^q$. As before, let 
$\phi\colon G\surj H$ be an epimorphism. 

\begin{prop}[\cite{CS}, \cite{Tr}]
\label{prop1}
Let  $k$ be a positive integer.  Then $\E_{k-1}(B _{\phi}) \cdot I^{q-1}
\subset \E_k(A _{\phi})$. Moreover, there is a positive integer $r$ such 
that $\E_k(A _{\phi})\cdot I^{r+k-q} \subset \E_{k-1}(B _{\phi})$.
\end{prop}
 
\begin{corollary}
\label{cor:cst}
For any integer $k\ge 1$, the subschemes of  
$\T_H$ given by  $V(\E_{k-1}(B _{\phi}))$ and $V(\E_{k}(A _{\phi}))$ 
coincide away from the identity element $1 \in \T_H$.
\end{corollary}

The most important case for us is the projection 
onto the maximal torsion-free abelian quotient, 
$\phi_{\abf}\colon G \to G_{\abf}$,   where 
$G_{\abf}=G_{\ab}/{\Tors}$.  The resulting covering is 
the maximal torsion-free abelian cover of $X$, and the 
corresponding Alexander modules are the classical objects 
considered in the theory of knots and links; see e.g.~\cite{Hi}, 
\cite{MM}. In the sequel, whenever $\phi$ is not mentioned, we will 
assume that $\phi$ is the  morphism $\phi_{\abf}$, and will 
denote the corresponding $\Z{G_{\abf}}$-modules simply 
by $B_G$ and $A_G$.

\subsection{Characteristic varieties}
\label{ss13}
As before, let $X$ be a connected CW-complex with finite $1$-skeleton, 
and $G=\pi_1(X,x_0)$. Consider the character group of $G$, 
i.e., the complex algebraic group 
\[
\T_G= \Hom (G, \C^*)= \Hom (G_{\ab}, \C^*).
\]  
Its identity component, $\T_G^0$, 
may be identified with the complex torus 
$\Hom (G_{\abf}, \C^*)= (\C^*)^n$, where $n=b_1(G)$. 
Each character $\rho\colon G\to \C^*$ gives rise to a 
rank~$1$ local system on $X$, denoted by $\C_{\rho}$. 
The {\em characteristic varieties} of $X$ are the 
jumping loci for homology with coefficients 
in such local systems:
\begin{equation*}
\label{eq:defchar}
\V_k (X)= \{ \rho\in \T_G  \mid \dim_{\C} H_1(X; \C_{\rho})\ge k \}.
\end{equation*}
It is readily seen that these sets form a decreasing sequence of 
algebraic subvarieties of $\T_G$, depending only on the isomorphism 
type of the group $G$.  In particular, we may write $\V_k (G):=\V_k (X)$. 

The Lemma below is well-known---see for example Hironaka \cite{Hir}. 
For the reader's convenience, we include a short proof, in a slightly 
more general context. 

\begin{prop}
\label{prop:echar}
Let $G$ be a finitely generated group. Then, for all $k\ge 1$, 
\[
\V_k(G)\cap \T_G^0= V_{k+1}(A_G\otimes \C) 
\quad \text{and}\quad \V_k(G) = V_{k+1}(A_{\ab}\otimes \C)
\]
away from $1$, with equality at $1$ for $k<b_1(G)$.
\end{prop}
\begin{proof}
Let $\rho\in \T_G^0 \setminus \{ 1\}$ be a 
non-trivial character. By the discussion in \S\ref{ss:tprel}, we have 
\[
\dim_{\C} H_1(X; \C_{\rho})\ge k 
\twoarr
\dim_{\C} (\C_{\rho}\otimes_{\C G_{\abf}} A_G\otimes \C)\ge k+1. 
\]
This is turn is equivalent to
$\rho\in V_{k+1}(A_G\otimes \C)$, by the definitions from 
\S\S\ref{subsec:fitt} and \ref{subsec:char var}. 

It is readily checked that the identity belongs to $\V_k(G)$ 
and to $V_{k+1}(A_G\otimes \C)$, whenever $k<b_1(G)$. 
The proof for $A_{\ab}$ is similar.
\end{proof}

As a corollary, we obtain the following result (see also Delzant \cite{De} 
for another proof). 

\begin{corollary} 
\label{cor:metabelian}
The characteristic varieties $\V_k(G)$ depend only on $G/G''$, the 
maximal metabelian quotient of $G$. 
\end{corollary}

\begin{proof}
Clearly, the projection $G\surj G/G''$ induces the identity on abelianizations. 
By \eqref{eq:bab}, the $\Z{G}_{\ab}$-module $B_{\ab}=G'/G''$ is 
determined by $G/G''$. By naturality of the Crowell exact 
sequence \eqref{eq:crow2}, the Alexander module $A_{\ab}$ also 
depends only on $G/G''$.  The conclusion follows from 
Proposition \ref{prop:echar}.   
\end{proof}

Let $H$ be a finitely-generated, free abelian group, and $\phi\colon G\surj H$ 
an epimorphism. 
The induced map, $\phi^*\colon \T_H \hookrightarrow \T_G$, is an inclusion,  
mapping $\T_H$ to a subtorus of  $\T_G^0$. From Corollary \ref{cor:cst} and 
Proposition \ref{prop:echar}, we obtain the following. 

\begin{corollary} 
\label{cor:jumpb}
Let $\rho$ be a non-trivial character in $\T_H$.  For all $k\ge 1$, 
\[
\phi^*(\rho) \in \V_k(G) \twoarr
\rho\in V(\E_{k-1}(B_{\phi})).
\]
\end{corollary}

\subsection{Finitely presented groups}
\label{subsec:fp}
Recall $X$ is a connected CW-complex with finite $1$-skeleton, 
and so $G=\pi_1(X)$ is finitely generated. For our purposes, we 
may assume without loss of generality that $G$ is, in fact, 
finitely presented.  

Indeed, we may  
construct a finitely presented group $\bar{G}$, together with an 
epimorphism, $f\colon \bar{G}\surj G$, so that the induced map 
on abelianizations, $f_*\colon \bar{G}_{\ab} \to G_{\ab}$, is an 
isomorphism. Given a commutative ring $S$, a ring homomorphism 
$\sigma\colon \Z{G}\to S$ puts a $\Z{G}$-module structure on $S$;   
moreover, the homology groups of $G$ with coefficients 
in this $\Z{G}$-module inherit an $S$-module structure.  
It is readily checked that the map
$\Hom(f,\id_S)\colon \Hom_{\rm Rings} (\Z{G}, S)\to 
\Hom_{\rm Rings} (\Z{\bar{G}}, S)$ is an isomorphism. 
Thus, $f_*\colon H_1(\bar{G}; S)\to H_1(G; S)$ is an 
$S$-isomor\-phism.  It follows that the Alexander modules 
and the characteristic varieties of $G$ and $\bar{G}$ are isomorphic. 

Now suppose $G= \langle x_1, \dots, x_m \mid r_1, \dots, r_h \rangle$ 
is a finite presentation.  Let $\partial r_i/\partial x_j $ 
be the Fox derivatives of the relators, viewed as elements in $\Z{G}$, 
and let  $\tilde{\phi}\colon \Z{G} \to \Z{H}$ be the linear extension of 
the homomorphism $\phi\colon G\to H$.  As shown by 
R.~Fox \cite{Fo}, the resulting matrix, 
\begin{equation*}
\label{eq:alex mat}
\AA_G^{\phi}=\big( \tilde\phi ( \partial r_i/\partial x_j )\big)\colon 
\Z{H}^h \to \Z{H}^m,
\end{equation*}
is a presentation matrix for the Alexander module $A_{\phi}$.  
In view of Proposition~\ref{prop:echar}, this yields an algorithm for 
computing the characteristic varieties $\V_k(G)$, as the 
determinantal varieties of the matrix $\AA_G^{\ab}$. The matrix 
$\AA_G=\AA_G^{\abf}$ will be called the {\em Alexander matrix} of $G$. 

We shall have occasion to consider various base 
changes on the  matrix $(\partial r_i/\partial x_j)$.  
If $\psi\colon \Z{G}\to S$ 
is a ring homomorphism, then applying $\psi$ to each entry 
yields a new matrix, $\AA^{\psi}_G \colon S^h \to S^m$. 

\section{Codimension-one components of characteristic varieties}
\label{sect:alexvc1}

In this section, we analyze the extent to which the Alexander 
polynomial of a space $X$ determines the codimension~$1$ 
irreducible components of the characteristic variety $\V_1(X)$. 

\subsection{The codimension-one stratum of the support variety}
\label{ss21}
Let $\E$ be a non-zero ideal in 
$\Lambda_q\otimes \C= \C[t_1^{\pm 1}, \dots, t_q^{\pm 1}]$, and 
set $\Delta = \gcd (\E)$. Denote by $V(\Delta)$ and $V(\E)$ the 
corresponding zero sets in $(\C^*)^q$. Let $\w(\E)$ be the 
union of all codimension-one irreducible components of $V(\E)$.

\begin{lemma}
\label{lem:edelta}
With the above notation, $\w(\E)=V(\Delta)$.
\end{lemma}

\begin{proof}
We use the fact that $\Lambda_q$ is a unique factorization domain. 
Let $W$ be a codimension-$1$ component of $V(\E)$. Then $W=V(g)$, 
for some irreducible Laurent polynomial $g\in \Lambda_q$. Since 
$V(g)\subset V(\E)$, the polynomial $g$ must divide $\Delta$, and 
so $W\subset V(\Delta)$.  

Conversely, note that $\Delta\ne 0$, since $\E\ne 0$. If $\Delta$ is a unit, 
then $V(\Delta)=\emptyset$. Otherwise, let $\{ f_j\}$ be the prime factors 
of $\Delta$, and $V(\Delta)=\bigcup_j V(f_j)$ the irreducible decomposition. 
It follows that $V(f_j)\subset V(\E)$, for all $j$. Since $\dim V(\E)\le n-1$, 
we infer that $V(\Delta)\subset \w(\E)$, and we are done.
\end{proof}

\subsection{Alexander polynomial} 
\label{subsec:alex poly}
Let $X$ be a CW-complex with finite $1$-skeleton, and let 
$G=\pi_1(X,x_0)$ be its fundamental group. Suppose we 
are given a surjective homomorphism, $\phi\colon G\surj H$, 
from $G$ to a torsion-free, finitely generated abelian group $H$. 
Define the {\em Alexander polynomial} of $X$ (relative to $\phi$) 
as
\[
\Delta_{\phi}^X:= \Delta_1(A_{\phi})\in \Z{H}. 
\]
From the discussion in \S\ref{sect:prel}, it follows that 
$\Delta_{\phi}^X$ depends only on $G$, modulo units in 
$\Z{H}$; consequently, we will write the Alexander polynomial 
as $\Delta_{\phi}^G$. In the case when $H=G_{\abf}$ and 
$\phi=\phi_{\abf}$, we will simply write the Alexander polynomial 
as $\Delta^G$, and view it as an element in $\Lambda_n$, 
where $n=b_1(G)$. 

Denote by $V(\Delta^G)$  the hypersurface of $\T_G^0=(\C^*)^n$ 
defined by the Alexander polynomial of $G$. Furthermore, denote 
by $\vv (G)$ the union of all codimension-one irreducible components 
of $\V_1(G)$ contained in $\T_G^0$.

Using Proposition \ref{prop:echar}, Lemma \ref{lem:edelta}, and the 
discussion at the end of \S\ref{subsec:fitt}, we obtain the following 
relationships between $\vv (G)$ and $\Delta^G$. 

\begin{corollary}
\label{cor:deltacd1}
Let $G$ be a finitely generated group.  
\begin{enumerate}
\item \label{dc1}
$\Delta^G=0 \twoarr \T_G^0\subset \V_1(G)$.
In this case, $\vv (G)=\emptyset$.
\item \label{dc2}
If $b_1(G)\ge 1$ and $\Delta^G\ne 0$, then
\begin{equation*}
\vv (G) =\begin{cases}
V(\Delta^G) & \text{if $b_1(G)>1$}\\
 V(\Delta^G)\coprod \{ 1\}  & \text{if $b_1(G)=1$}.
\end{cases}
\end{equation*}
\item \label{dc3} If $b_1(G)\ge 2$, then 
\[
\vv (G)=\emptyset  \twoarr \Delta^G\doteq \const. 
\]
\end{enumerate}
\end{corollary}

\subsection{The Alexander polynomial of a hyperplane arrangement}
\label{subsec:hyp arr}

Here is a quick application of this result. 
Let $\A$ be a hyperplane arrangement in 
$\C^{\ell}$, with complement $X(\A)=\C^{\ell}\setminus 
\bigcup_{H\in A} H$.  If $\A'$ is a generic two-dimensional 
section of $\A$, a well-known theorem 
guarantees that the inclusion $X(\A')\inj X(\A)$ induces an 
isomorphism $\pi_1(X(\A'))\isom \pi_1(X(\A))$.    
So, for the purpose of studying characteristic varieties 
and Alexander polynomials of arrangement groups, we 
may as well work with (affine) line arrangements.  

\begin{example}
\label{ex:pencil}
Let $\A$ be a pencil of $n$ lines through the origin of $\C^2$, 
defined by the equation $x^n-y^n=0$. The fundamental group 
of the complement, $G=\pi_1(X(\A))$, is isomorphic to 
$\Z\times F_{n-1}$, where $F_m$ denotes the free group on $m$ generators. 
If $n=1$ or $2$, then 
$\Delta^G=1$ and $\V_1(G)=\{1\}$, but otherwise,  
$\Delta^G= (t_1\cdots t_n-1)^{n-2}$ and 
$\V_1(G)=\vv(G)=\{ t_1\cdots t_n-1=0\}$.
\end{example}

This example is very special, though: in general, arrangement 
groups have trivial Alexander polynomial. 

\begin{prop}
\label{prop:delta arr}
Let $\A$ be an arrangement of $n$ lines.  The Alexander polynomial 
of $G=\pi_1(X(\A))$ is non-constant if and only if $\A$ is a 
pencil and $n\ge 3$.
\end{prop}

\begin{proof}
Suppose $\A$ is not a pencil.  
By the above, each sub-arrangement which is a pencil 
of $k\ge 3$ lines determines a ``local'' component of 
$\V_1(G)$ of dimension $k-1$.  If $n\le 5$, then 
it is readily seen that all components of $\V_1(G)$ 
must be local. In general, there are other irreducible 
components, but each non-local component 
has dimension at most $4$, by \cite[Theorem 7.2]{PY}.
Since $\A$ is not a pencil, all components of $\V_1(G)$ 
must have codimension at least $2$, i.e., 
$\vv(G)=\emptyset$. Since $b_1(G)\ge 2$,  
Corollary \ref{cor:deltacd1}\eqref{dc3} implies 
$\Delta^G\doteq \const$. 
\end{proof}

\subsection{Parallel subtori}
\label{susbec:parallel}
Let $H$ be a finitely generated abelian group.  
By a {\em subtorus} of $\T_H$ we mean a subvariety $W$ 
of the form 
\[
W=\rho \cdot S,
\] 
where $\rho\in \T_H$ and $S=f^*(\T_D)$, for some epimorphism
$f\colon H\surj D$ to a free abelian group $D$. Clearly, the 
translation factor $\rho$ is not uniquely determined by $W$, 
while the subtorus $S$ is.  We call this subtorus the {\em direction} 
of $W$, and write $S=\dire (W)$. 

Two subtori, $W$, $W' \subset \T_H$, are said to be {\em parallel} 
if the direction of one is contained in the direction of the other, i.e., 
\[
W \parallel W' \twoarr 
\dire (W)\subset \dire(W'), \text{ or }
\dire (W')\subset \dire(W).
\]
Clearly, two codimension-one subtori are 
parallel if and only if they have the same direction.

\subsection{Single essential variable}
\label{subsec:ess var}
Assume $q:=\rank H$ is 
positive.  In a suitable coordinate system $(t_1,\dots ,t_q)$ 
on $\T_H^0=(\C^*)^q$, a codimension-one subtorus contained 
in $\T_H^0$ is given by an equation of the form $t_1-z=0$, 
for some $z\in \C^*$.  This observation motivates the following 
definition. 

\begin{definition}
\label{def:essvar}
Let $H$ be a finitely generated free abelian group.
We say an element $\Delta\in \Z{H}$ has {\em a single essential 
variable} if 
\[
\Delta \doteq \tilde\nu(P),
\]
for some Laurent polynomial $P\in \Z[\Z]$, and some 
homomorphism $\nu \colon \Z\to H$. 
\end{definition}

Now assume $q:=\rank H$ is positive.
Note that, in the above definition, $\nu$ may be supposed to  
be a split monomorphism. Identifying $\Z H$ with the ring of
Laurent polynomials, 
$\Lambda_q=\Z[t_1^{\pm 1}, \dots , t_q^{\pm 1}]$,
Definition \ref{def:essvar} can be understood in more 
concrete terms, as follows.

\begin{lemma}
\label{lem:newton}
For a Laurent polynomial $\Delta\in \Lambda_q$, the following 
conditions are equivalent.
\begin{enumerate}
\item \label{np1}
$\Delta$ has a single essential variable.
\item \label{npc}
There is a Laurent polynomial $P\in \C[\Z]$ and a group
homomorphism $\nu \colon \Z\to H$, such that 
$\Delta \doteq_{\C} \tilde\nu(P)$.
\item  \label{np2}
$\Delta (t_1,\dots , t_q)\doteq P(t_1^{e_1}\cdots t_q^{e_q})$, 
for some polynomial $P\in \Z[t^{\pm 1}]$, and some coprime 
exponents $e_1,\dots, e_q\in \N$. 
\item  \label{np3}
The Newton polytope of $\Delta$ is a line segment. 
\end{enumerate}
\end{lemma}

The proof of this Lemma is left as an exercise in the definitions. 
We are now ready to connect this `$1$-dimensionality' property 
of the Alexander polynomial to the geometry of the codimension-one 
stratum of the characteristic variety. 

\begin{prop}
\label{prop:cod1ess}
Let $G$ be a finitely generated group, with $b_1(G)\ge 2$. 
The Alexander polynomial $\Delta^G$ has a single essential 
variable if and only if either $\vv (G)=\emptyset$, or all the irreducible 
components of $\vv (G)$ are parallel, codimension-one subtori of 
$\T_G^0$.
\end{prop}

\begin{proof}
By Corollary \ref{cor:deltacd1}\eqref{dc3}, we may assume  that 
$\Delta^G \not\doteq \const$ and $\vv (G)\neq \emptyset$. 

Set $n=b_1(G)$, and suppose $\Delta^G\in \Lambda_n$ has a 
single essential variable. Then 
$\Delta^G \doteq_{\C} \prod_j (t_1-z_j)^{m_j}$, for some pairwise 
distinct $z_j\in \C^*$. By Corollary \ref{cor:deltacd1}\eqref{dc2}, 
\begin{equation}
\label{eq:pdec}
\vv (G)= \bigcup_j \{ t_1=z_j \}\, ,
\end{equation}
as asserted. 

Conversely, let $\Delta^G \doteq_{\C} \prod_k f_k^{\mu_k}$ be the 
decomposition of $\Delta^G$ into irreducible factors, leading to the 
decomposition into irreducible components for its zero locus,
\begin{equation}
\label{eq:alexdec}
V(\Delta^G)= \bigcup_k \{ f_k= 0 \}. 
\end{equation} 
By assumption, \eqref{eq:pdec} also holds, in a suitable coordinate system 
$(t_1,\dots ,t_n)$ on $\T_G^0=(\C^*)^n$. 
Comparing decompositions \eqref{eq:pdec} and \eqref{eq:alexdec}, 
and making use of Corollary \ref{cor:deltacd1}\eqref{dc2}, we infer that 
$\Delta^G \doteq_{\C} \prod_j (t_1-z_j)^{\mu_j}$. This finishes 
the proof.
\end{proof}

\section{Quasi-projective groups}
\label{sect:qprojalex}

In  this section, we analyze the characteristic variety $\V_1(G)$ and 
the Alexander polynomial $\Delta^G$ 
in the case when $G$ can be realized as the fundamental group 
of a smooth, connected, quasi-projective complex variety.

\subsection{Admissible maps and isotropic subspaces}
\label{subsec:qp}
We build on a foundational result of Arapura \cite[Theorem V.1.6]{A} 
on the structure of the first characteristic variety of a quasi-projective 
group, together with some refinements from \cite{D06} and \cite{DPS}.

\begin{theorem}[Arapura \cite{A}]
\label{thm:arapura}
Let $G= \pi_1(M)$ be a quasi-projective group.  Then all irreducible components 
$W$ of $\V_1(G)$ are subtori. More precisely, if $\dim W>0$, then 
\begin{equation*}
\label{eq:arapura}
W=\rho \cdot f^*\big(\T_{\pi_1(C)}\big),
\end{equation*}
where $\rho\in \T_G$ is a torsion element, and $f\colon M\to C$ is 
a surjective regular map with connected generic fiber onto a smooth 
complex curve $C$.
\end{theorem}

A map $f\colon M\to C$ as above is called {\em admissible}.  
Let $F$ be the generic fiber of $f$, and $i\colon F\to M$ the inclusion.  
The induced homomorphism, $f_{\sharp}\colon \pi_1(M) \to \pi_1(C)$,  
is surjective.  Furthermore, the composite 
$\xymatrixcolsep{16pt}
\xymatrix{\pi_1(F) \ar[r]^{i_{\sharp}} & \pi_1(M) \ar[r]^(.55){\rho} &\C^*}$ 
is trivial. On the other hand, the sequence 
$\xymatrixcolsep{16pt} 
\xymatrix{H_1(F) \ar[r]^{i_*} & H_1(M) \ar[r]^{f_*} &H_1(C) \ar[r]& 0}$ 
is not always exact in the middle. As shown in \cite{D06}, the quotient
\begin{equation}
\label{eq:tf}
T(f)=\ker(f_*)/\im(i_*)
\end{equation}
is a finite abelian group.  Furthermore, the components of $\V_1(G)$ with  
direction $f^*(\T_{\pi_1(C)})$ are parametrized by the 
Pontryagin dual $\widehat{T}(f)=\Hom(T(f),\C^*)$, except when 
$\chi(C)=0$ (i.e., $C$ is either an elliptic curve, or $C=\C^*$), in 
which case the trivial character must be excluded from 
$\widehat{T}(f)$. 

We also need another definition, taken from \cite{DPS}. Let  
\begin{equation}
\label{eq:cup prod}
\cup_M \colon H^1(M;\C)\wedge H^1(M;\C) \to H^2(M;\C),
\end{equation}
be the cup-product map.  Given a subspace 
$H\subset H^1(M;\C)$, denote by $\cup_H$ the restriction of 
$\cup_M$ to $H\wedge H$. We say that $H$ is  $0$-{\em isotropic} 
if $\cup_H=0$, and $1$-{\em isotropic} if $\dim \im (\cup_H)=1$ 
and $\cup_H$ is a non-degenerate linear skew-form on $H$. 
Clearly, both notions depend only on $G=\pi_1(M)$. Note also 
that $H^1(M; \C)$ is naturally identified with the tangent space 
to $\T_G$ at the identity, $T_1\T_G$.

\subsection{Position obstructions}
\label{subsec:pos obs}
Our next result gives obstructions on the position of all components 
of $\V_1(G)$, thus extending Theorem B from \cite{DPS}, where only 
components through $1$ were considered. Part \eqref{a1} concerns 
the ambient position of each component in $\T_G$, while 
Parts \eqref{a2}--\eqref{a3} describe the relative position of 
each pair of components.

\begin{theorem}
\label{thm:posobs}
For an arbitrary quasi-projective group $G$, the following hold.
\begin{enumerate}

\item \label{a1}
If $W$ is an irreducible component of $\V_1(G)$, then 
$W\subset \T_G $  is a subtorus, and the subspace 
$T_1 \dire (W)\subset T_1\T_G$ is either $0$-isotropic 
or $1$-isotropic.

\item \label{a2}
If $W$ and $W'$ are two distinct components of 
$\V_1(G)$, then either $\dire (W) = \dire (W')$, or 
$T_1 \dire (W) \cap T_1 \dire (W')= \{ 0\}$. 

\item \label{a3}
For each pair of distinct components, $W$ and $W'$, the 
intersection $W\cap W'$ is a finite (possibly empty) set.
\end{enumerate}
\end{theorem}

\begin{proof}
Plainly, we may assume $b_1(G)\ge 1$, and consider only 
positive-dimensional components.   

Part \eqref{a1}.  We know from Theorem~\ref{thm:arapura} that 
$W$ is a subtorus of the form $W=\rho \cdot S$, with $\rho$ 
a finite-order character, and $S:= \dire (W)=f^*(\T_{\pi_1(C)})$, 
for some admissible map $f\colon M\to C$. 
If $C$ is compact of positive genus and 
$f^*\colon H^2(C;\C)\to H^2(M;\C)$ is a 
monomorphism, then the subspace $T_1 S=f^*(H^1(C;\C))$ 
is $1$-isotropic; otherwise, this subspace is $0$-isotropic.

Part \eqref{a2}. As above, write $W'=\rho' \cdot S'$, 
with $S'= \dire (W')=f'^*(\T_{\pi_1(C')})$, for some admissible 
map $f'\colon M\to C'$.  If $W$ and $W'$ are not parallel, 
then the arguments from \cite[Lemmas 6.3--6.4]{DPS} 
imply that $S \cap S'$ must be finite, and so 
$T_1 S\cap T_1 S'=\{0\}$. 

So suppose $W$ and $W'$ are parallel, but $S\ne S'$. 
We then have a strict inclusion, say $S\subset S'$. 
This forces $\dim S' \ge 2$, and so $b_1(C')=\dim W' \ge 2$.  
Hence, we can apply \cite[Corollary 4.6]{D06} to the admissible 
map $f'\colon M\to C'$ and the rank one local system 
$\LL_1=\C_{\rho}$.\footnote{In \cite{D06},  
``for a generic local system" actually means ``for all but 
finitely many local systems".} 
Set $W_1=\rho \cdot S'$. Since $W \subseteq W_1$, the 
restriction of $\LL_1$ to the generic fiber of $f'$ must be trivial. 
Thus, $\LL_1$ corresponds to a character
$\alpha\in \widehat{T}(f')$. There are two cases to consider.

{\em Case 1. Either $C'$ is not an elliptic curve, 
or $C'$ is an elliptic curve and the character $\alpha$ 
is non-trivial.} 
Then, according to \cite[Corollary 5.8]{D06}, $W_1$ is an irreducible 
component of $\V_1(G)$. Since  $W$ is also an irreducible component 
of $\V_1(G)$, it follows that $W=W_1$, and hence $S=S'$, a contradiction.

{\em Case 2. $C'$ is an elliptic curve and 
the character $\alpha$ is trivial.} 
Then, again by \cite[Corollary 5.8]{D06}, $W_1=S'$ is not an 
irreducible component of $\V_1(G)$.  In this case $\dim S'=2$ 
and the elements in $T_1 S'=f'^*(H^1(C';\C))$ have Hodge 
type $(1,0)$ and $(0,1)$.  On the other hand, $\dim S=1$, 
which implies that $C=\C^*$, and the elements in 
$T_1 S=f^*(H^1(\C^*;\C))$ have Hodge type $(1,1)$.  
This contradicts $T_1 S \subset T_1S'$.

Part \eqref{a3}. This is an immediate consequence of Part \eqref{a2}.
\end{proof}

Note that, in Part~\eqref{a1}, both $0$-isotropic and $1$-isotropic 
subspaces may arise; see \cite{DPS}.  The second situation 
described in Part~\eqref{a2} often occurs for complements 
of complex hyperplane arrangements, see for instance \cite{S}. 
The first situation often occurs for complements of Seifert links, 
see for instance Example \ref{ex4}.

\subsection{Alexander polynomial and quasi-projectivity}
\label{subsec:alexone}
The relative position obstruction from Theorem \ref{thm:posobs}\eqref{a2} 
translates into restrictions on the multivariable Alexander polynomial of 
a quasi-projective group, as follows. 

\begin{theorem}
\label{thm:alexone}
Let $G=\pi_1(M)$ be the fundamental group of a smooth, connected, 
complex quasi-projective variety. 
\begin{enumerate}
\item \label{b1}
If $b_1(G)\ne 2$, then the Alexander polynomial $\Delta^G$ has a 
single essential variable. 
\item \label{b2}
If $b_1(G)\ge 2$, and $\Delta^G$ has a single essential variable, 
then either 
\begin{enumerate}
\item $\Delta^G = 0$, or 
\item $\Delta^G(t_1,\dots , t_n) \doteq_{\C} P(u)$, where $P$ is 
a product of cyclotomic polynomials (possibly equal to $1$), 
and $u=t_1^{e_1}\cdots t_n^{e_n}$, with $\gcd(e_1,\dots,e_n)=1$. 
\end{enumerate}
\item \label{b3}  If $M$ is actually a projective variety, then 
$\Delta^G \doteq \const$. 
\end{enumerate}
\end{theorem}

\begin{proof}
In Part \eqref{b1}, we may assume $b_1(G)\ge 3$. The 
claim then follows from Proposition \ref{prop:cod1ess} and 
Theorem \ref{thm:posobs}\eqref{a2}.

In Part \eqref{b2}, we may assume---after a  change of 
variables if necessary---that $\Delta^G \doteq P(t_1)$, 
with $P$ a non-constant polynomial in $\Z[t_1]$. By 
Corollary \ref{cor:deltacd1}\eqref{dc2} and Arapura's 
Theorem \ref{thm:arapura}, all roots of $P$ are roots of unity.

In Part \eqref{b3}, since the odd Betti numbers of $M$ are even, 
we may assume $b_1(G)\ge 2$, for, otherwise, clearly $\Delta^G$ is constant. 
By Corollary \ref{cor:deltacd1}\eqref{dc3}, we are left with 
showing that $\vv (G)=\emptyset$.  Suppose, to the contrary, 
that $\V_1(G)$ has a codimension-one component, call it $W$;  
then there is an admissible map, $f\colon M\to C$, onto a 
projective curve, such that $\dire (W)= f^*\big(\T_{\pi_1(C)}\big)$. 
This implies $b_1(C)$ is odd, a contradiction.
\end{proof}

\begin{corollary}
\label{cor:alexshort}
If $G$ is a quasi-projective group with $n=b_1(G)\ge 3$, then 
\[
\Delta^G (t_1,\dots , t_n)\doteq c P(t_1^{e_1}\cdots t_n^{e_n}),
\] 
where $P\in \Z[t]$ is a product of cyclotomic polynomials, 
$c\in \Z$, and $\gcd(e_1,\dots ,e_n)=1$. 
\end{corollary}

\begin{remark}
\label{rem:lib}
A result similar to Theorem \ref{thm:alexone}\eqref{b2} was 
obtained by Libgober in \cite{Li1} for complements of plane 
algebraic curves, $M=\PP^2\setminus C$.  Yet his result only 
holds for the {\em single variable} Alexander polynomial, 
$\Delta=\Delta^G_{\phi_{\lk}}$, associated to $G=\pi_1(M)$ and  
the total linking homomorphism, $\phi_{\lk}\colon G\to \Z$. 
Note that the polynomial $\Delta(t) \in \Z[t^{\pm 1}]$ carries 
information only on the intersection of a one-dimensional 
subtorus of $\T_G$ with $\V_1(G)$, whereas our general 
result from Theorem \ref{thm:alexone}, when combined with 
Proposition \ref{prop:cod1ess}, puts severe restrictions on 
the global structure of $\V_1(G)$.  
\end{remark}

\subsection{The case $b_1=2$} 
\label{subsec:betti 2}
Concerning Part \eqref{b1} of Theorem \ref{thm:alexone}, we 
now give an example showing that the Alexander polynomials of 
quasi-projective groups with $b_1=2$ do exhibit an exceptional 
behavior.

\begin{example}
\label{ex:bdoi}
Consider the surface $X\colon xy-z^2=0$ in $\C^3$. Blowing up the 
origin, we obtain a smooth surface $Y$, which is known to be the 
total space of the line bundle $\OO(-2)$ on $\PP^1$. Hence
$\tilde{H}_0(Y;\Z)=H_1(Y;\Z)=0$.
Indeed, the smooth projective conic $X'\colon xy-z^2=0$ in 
$\PP^2$ is rational of degree $2$ (use the genus-degree formula).
Consider the lines $L_1\colon x=1$ and $L_2\colon y=1$ on $X$ 
and denote by $C_1$ and $C_2$ their proper transforms
in $Y$. Since the lines $L_j$ do not pass through the origin, 
the curves $C_j$ are isomorphic to $L_j$, therefore to $\C$.

Let $M=Y \setminus (C_1 \cup C_2)$.  The above discussion 
implies that $H_1(M;\Z)=\Z^2$, with basis given by some small 
loops $\gamma_1$ and $\gamma_2$ about $C_1$ and $C_2$.  
Consider the group $G=\pi_1(M)$.  
Clearly, $G$ is quasi-projective, and $b_1(G)=2$. We claim 
that the Alexander polynomial of $G$ does {\em not} have a 
single essential variable. 

Indeed, let $C=\C \setminus \{1\}$, and consider the admissible 
map $f\colon M \to C$, obtained by composing the projection map 
$X \setminus (L_1 \cup L_2) \to C$, $(x,y,z)\mapsto x$ with the 
blow-up morphism.  The induced homomorphism, 
$f^*\colon H^1(C;\Z) \to H^1(M;\Z)$, 
has image equal to $\Z \cdot \gamma_1^*$.
Recall from \S\ref{subsec:qp} that $T(f)$ denotes the finite abelian 
group $\ker (f_*)/\im (i_*)$, where $i\colon F\hookrightarrow M$ 
is the inclusion of the generic fiber of $f$. In our situation, 
a computation shows that $T(f)=\Z/2\Z$.
It follows from \cite{D06} that there is a unique 
$1$-dimensional component of $\V_1(G)$ associated to $f$, 
call it $W_f$, with tangent direction $\gamma_1^*$. 
Proceeding in exactly the same way, starting from the 
other projection map, $g\colon (x,y,z)\mapsto y$, we 
obtain another $1$-dimensional component, call it 
$W_g$, with tangent direction $\gamma_2^*$.
Clearly, $W_f\nparallel W_g$. By Proposition \ref{prop:cod1ess}, 
$\Delta^G$ does not have a single essential variable.
\end{example}

\subsection{Boundary manifolds of line arrangements}
\label{subsec:bdry mfd}
As an application of our methods, we now give a class of 
examples where the restrictions imposed by 
Theorem \ref{thm:alexone} completely settle Serre's problem. 
More applications will be given in a forthcoming paper.
 
Let $\A=\{\ell_0,\dots ,\ell_n\}$ be an arrangement of 
lines in $\PP^2$.  Of course, the complement 
$X(\A)=\PP^2\setminus \bigcup_{i=0}^{n} \ell_i$ is a smooth, 
quasi-projective variety, and so $\pi_1(X(\A))$ is a quasi-projective 
group. On the other hand, consider the closed, orientable 
$3$-manifold $M(\A)$, defined as the boundary of a regular 
neighborhood of the curve $C=\bigcup_{i=0}^{n} \ell_i $ in 
$\PP^2$. The fundamental group of the boundary manifold, 
$G_{\A}= \pi_1(M_{\A})$, admits a commutator-relators 
presentation, with $b_1(G_{\A})\ge n$, see \cite{CS06}. 

\begin{prop}
\label{prop:serrebd}
The group $G_{\A}$ is quasi-projective if and only if $\A$ is a pencil 
or a near-pencil.
\end{prop}

\begin{proof}
If $\A$ is a pencil of $n+1$ lines, then $G_{\A}$ is the free group 
of rank $n$; if $\A$ is a near-pencil (i.e., a pencil of $n$ lines, 
together with an extra line in general position), then $G_{\A}$ is 
the product of $\Z$ with the fundamental group of a compact 
Riemann surface of genus $n-1$, see \cite{CS06}.   In either 
case, $G_{\A}$ is quasi-projective. 

For the converse, let $\A$ be an arrangement that is not a 
pencil or a near-pencil; in particular, $n\ge 3$. By 
\cite[Theorem 5.2]{CS06}, the Alexander polynomial 
$\Delta^{G_{\A}}$ is a product of (irreducible) factors of the 
form $t_{i_1}\cdots t_{i_r}-1$; by \cite[Proposition 5.5]{CS06}, 
$\Delta^{G_{\A}}$ must have at least 
two distinct such factors.

Now suppose $G_{\A}$ is quasi-projective.  Since 
$b_1(G_{\A})\ge 3$, we infer from 
Theorem \ref{thm:alexone}\eqref{b1} that 
$\Delta^{G_{\A}}$ must have a single essential variable. 
By Proposition \ref{prop:cod1ess} and 
Corollary \ref{cor:deltacd1}\eqref{dc2}, the components of 
$V(\Delta^{G_{\A}})$ must be parallel, codimension-one 
subtori. Since there are at least two such subtori containing $1$, 
we reach a contradiction.
\end{proof}

\section{Multiplicities, twisted Betti ranks, and generic bounds}
\label{sect:multupper}

We saw in Section \ref{sect:alexvc1} that 
the {\em reduced} Alexander polynomial 
of $G$ determines the codimension-one stratum or the first 
characteristic variety, $\vv(G)$. In this section, we pursue this 
approach, analyzing the connection between the multiplicities 
of the factors of $\Delta^G$ and the higher-depth characteristic 
varieties, $\V_k(G)$.

\subsection{Generic Betti ranks}
\label{subsec:generic betti}

Let $G$ be a finitely generated group. Set $\Lambda:=\C[G_{\abf}]$. 
Consider the irreducible subvariety $V(\p) \subset \T_G^0$ 
associated to a prime ideal $\p\subset \Lambda$. 

\begin{definition}
\label{def:bgen}
The {\em generic Betti number} of $G$ relative to $\p$ is 
\[
\bg (G, \p)= \max  \{ k \mid  V(\p) \subset \V_k(G) \} .
\]
When $\p =(f)$, we abbreviate $\bg (G, (f))$ to $\bg (G, f)$.
\end{definition}

For a character $\rho\in \T_G$, set $b_1(G, \rho)
= \dim_{\C} H_1(G; \C_{\rho})$. We then have 
\begin{equation}
\label{eq:rhogen}
b_1(G, \rho)=\bg (G, \p) ,
\end{equation}
for $\rho$ in the Zariski open, non-empty subset
$V(\p)\setminus \V_{b+1}(G)$, where $b=\bg (G, \p)$. 
See also Farber \cite[Theorem 1.5]{Fa} for another interpretation. 
As the next example shows, generic Betti numbers are {\em not} 
always determined by the Alexander polynomial. 

\begin{example}
\label{ex:alexvsgen}
By a classical result of Lyndon (see \cite{Fo}), one may construct 
groups with quite complicated Alexander matrices. 
Denote by $F_m$ the free group on $m$ generators. Suppose  
$v_1, \dots ,v_m$ are elements in $\Z[\Z^m]$ satisfying 
$\sum_{j=1}^m (x_j -1) v_j=0$. There exists then an element 
$r\in F_m'$ such that $v_j= (\partial r/\partial x_j)^{\ab}$, for all $j$. 

Given any two elements $\Phi$, $\Psi\in \Z[\Z^3]$, we may 
construct in this way a commutator-relators group, 
$G=\langle x_1,x_2,x_3 \mid r_1, r_2\rangle$, 
with Alexander matrix
\[
\AA_G\begin{pmatrix}
(x_2-1)\Phi & (1-x_1)\Phi & 0\\ 
0 & (x_3-1)\Psi & (1-x_2)\Psi
\end{pmatrix}.  
\]
An easy computation reveals that $\Delta^G= (x_2-1)\Phi\Psi$. 

Now let $G_1$ and $G_2$ be groups corresponding to 
$\Phi_1=\Psi_1=\alpha \beta$, respectively 
$\Phi_2= \alpha^2$, $\Psi_2=\beta^2$, where 
$\alpha= x_1 x_2+1$ and $\beta= x_2 x_3+1$.  
These two groups share the same Alexander polynomial, 
$\Delta= (x_2-1)\alpha^2 \beta^2$. On the other hand, 
$\bg (G_1, \alpha)\ge 2$, whereas $\bg (G_2, \alpha)\le 1$.
\end{example}

\subsection{Generic upper bounds}
\label{subsec:generic ub}
Nevertheless, as the next result shows, the multiplicities 
of the Alexander polynomial  provide useful 
upper bounds for the generic Betti numbers.

First, some notation, which will be needed in the proof. Let 
$\OO_n$ be the ring of germs of analytic functions at $0\in \C^n$. 
The {\em order} of a germ $f\in \OO_n$ at $0$, denoted 
$\nu_0(f)$, is the degree of the first non-vanishing 
homogeneous term in the Taylor expansion of $f$ at $0$. 

\begin{theorem}
\label{thm:uppergen}
Let $\Delta^G \doteq_{\C} f_1^{\mu_1}\cdots f_s^{\mu_s}$ be the 
decomposition into irreducible factors of the Alexander polynomial
of a finitely generated group $G$. 
Then, for each $j$, 
\[
\bg (G, f_j)\le \mu_j .
\]
\end{theorem}

\begin{proof}
By the discussion from \S\ref{subsec:fp}, we may 
assume $G$ admits a finite presentation, say 
$G= \langle x_1, \dots, x_m \mid r_1, \dots, r_h \rangle$. 
Set $n=b_1(G)$, $\Delta= \Delta^G$, and identify $\Z{G}_{\abf}=\Lambda_n$. 

Fix an index $j\in \{1,\dots ,s\}$, and put $W=V(f_j)$. To prove 
the desired inequality, it is enough to find a Zariski open, 
non-empty subset $U\subset W$ such that 
$b_1(G, \rho)\le \mu_j$,  for every $\rho \in U$. 

By definition of $\Delta^G$, there exist 
$g_1,\dots , g_p\in \Lambda_n$ with $\gcd(g_1,\dots ,g_p)=1$ 
such that 
$\E_1(A_G)= (\Delta^G\cdot g_{1},\dots, \Delta^G\cdot g_{p})$. 
Define 
\begin{equation}
\label{eq:defu}
U= W\setminus \big( \{ 1\} \cup  \Sigma \cup Z \cup Y \big) ,
\end{equation}
where $Y=  V(g_{1})\cap \cdots \cap V(g_p)$, the set 
$\Sigma$ is the singular part of $W$, 
and $Z$ is the union of those irreducible components 
of $\V_1(G)$ contained in $\T_G^0$, yet different from $W$. 

Clearly, $U$ is a Zariski open set.  To show  
that $U$ is non-empty, we need to verify 
$W\not\subset  \{ 1\} \cup  \Sigma \cup Z \cup Y $. 
If $W\subset \{ 1\}$ then necessarily 
$n=1$ and $\Delta(1)=0$, which implies $n\ge 2$, 
a contradiction. Obviously, $W \not\subset \Sigma$. 
By Corollary \ref{cor:deltacd1}\eqref{dc2}, 
each $V(f_i)$ is a component of 
$\V_1(G)$ contained in $\T_G^0$; in particular,
$W \not\subset Z$. Finally, $W\subset Y$ would 
imply that $f_j$ divides each $g_{i}$, which is impossible.  

Now fix $\rho\in U$. Since $\rho\in W\setminus (\Sigma \cup Z)$, 
we may find analytic coordinates in a neighborhood $U_0\subset U$ 
of $\rho$, say $(z,w)$ with $z\in \OO_n$ and $w\in \OO_n^{n-1}$,  such that, in $U_0$,
\begin{equation*}
W= \{ z=0\},\quad
\{ w=0 \}\cap \V_1(G)\{\rho \},\quad
\big. \frac{d}{dz}\Big|_{z=0} f_j(z, 0)\ne 0.
\end{equation*}
In particular, $\nu_0 f_j(z, 0)=1$. 

Let  $\AA_G \colon \Lambda_n^h \to \Lambda_n^m$ be the 
Alexander matrix of $G$.  The base change 
$\tau\colon \Lambda_n \to \OO_1$, 
given by restriction of Laurent polynomials to $\{ w=0 \}$, 
yields the `local' Alexander matrix 
$\AA_G^{\tau}\colon \OO_1^h \to\OO_1^m$. 
For $t\in (\C, 0)$ belonging to the local transversal slice
to $W$ at $\rho$ given by $\{ w=0 \}$, 
since $\rho \ne 1$, we have
\begin{equation}
\label{eq:blocal}
1+  b_1(G, t) = \dim_{\C} \coker \AA_G^{\tau} (t),
\end{equation}
where $\AA_G^{\tau} (t) = \C_t\otimes_{\OO_1} \AA_G^{\tau}$ 
denotes the evaluation of the Alexander matrix at $t$.  
On the other hand, since $\OO_1$ is a PID, 
\begin{equation}
\label{eq:cokdec}
\coker \AA_G^{\tau}= \OO_1^r \oplus \Big(\bigoplus_{k=1}^\ell 
\OO_1/z^{\nu_k}\OO_1\Big) ,
\end{equation}
for some integers $r\ge 0$ and $\nu_k>0$.
By comparing \eqref{eq:blocal} and \eqref{eq:cokdec}, we infer 
that $r=1$ (by taking $t\ne 0$),  and $ b_1(G, 0)=\ell\ge 1$ 
(by taking $t=0$). 

Now, since $\rho\notin Y$, we must have 
$g_{i}(\rho)\ne 0$, for some $i$.  Thus 
$\tau (\Delta)\in \E_1(\coker \AA_G^{\tau})$, by base change. 
We deduce from formula \eqref{eq:delta pid}  
that $\tau(\Delta)\in z^{\nu_1+\dots +\nu_{\ell}} \OO_1$. 
In particular, $\nu_0 \Delta(z, 0)$ is at least $b_1(G, 0)=b_1(G, \rho)$.

Making use again of \eqref{eq:defu}, we see that $f_i(\rho)\ne 0$, 
for $i\ne j$; hence, $\nu_0 \Delta(z, 0)=\mu_j\cdot \nu_0 f_j(z, 0)$. 
But $\nu_0 f_j(z, 0)=1$, and so $\mu_j\ge b_1(G, \rho)$. 
This finishes the proof. 
\end{proof}

\subsection{Discussion of the generic bound}
\label{subsec:discuss genbound}

The inequality from Theorem \ref{thm:uppergen} may 
well be strict, as the following example shows. 

\begin{example}
\label{ex:gensharp}
Let $f\in \Z[\Z^2]$ be a Laurent polynomial, 
irreducible over $\C$. Using the method described in 
Example \ref{ex:alexvsgen}, we can construct, for each 
$k\ge 1$, a group $G_k=\langle x_1, x_2 \mid r_k\rangle$ 
with Alexander polynomial $\Delta^{G_k}=f^k$. Clearly, 
$\bg (G_k, f)=1$, while the multiplicity of $f$ as a factor 
of $\Delta^{G_k}$ is $\mu=k$. 
\end{example}

The next example shows that the inequality 
$\bg (G, f_j)\le \mu_j$ is sharp.

\begin{example}
\label{ex:central}
Let $G= \Z\times F_{n-1}$ be the fundamental group of the 
complement of a pencil of $n\ge 3$ lines in $\C^2$. 
Recall from Example \ref{ex:pencil} that $\Delta^G= f^{n-2}$, 
where $f=t_1\cdots t_n-1$.  It is readily checked that 
$\bg (G, f)=n-2$, which matches the multiplicity of $f$ 
as a factor of $\Delta^G$. 
\end{example}

Finally, let us note that the inequality from 
Theorem \ref{thm:uppergen} may well fail 
for {\em non-generic} twisted Betti ranks. 

\begin{example}
\label{ex:fail nongen}
Let $G=\langle x_1,x_2,x_3 \mid r_1, r_2\rangle$, be 
the group $G_1$ from Example \ref{ex:alexvsgen}, with 
Alexander polynomial $\Delta^G= (x_2-1)\Phi\Psi$. 
Consider the the character $\rho =(-1,1,-1)$, 
belonging to $V(x_2-1)\subset \vv (G)$. 
A computation shows that $b_1(G, \rho)=2$, whereas 
the multiplicity of the factor $x_2-1$ in $\Delta^G$ is $1$.
\end{example}

\section{Almost principal Alexander ideals and multiplicity bounds}
\label{sect:aai}

In this section, we delineate a class of groups 
$G$ (which includes the above example), for which the Alexander 
polynomial may be used to produce another upper bound for 
$b_1(G, \rho)$, valid for {\em all}\/ nontrivial local systems in 
$\T_G^0$.

\subsection{Almost principal ideals}
\label{subsec:aai}
We start with a definition, inspired by work of 
Eisenbud-Neumann \cite{EN} and McMullen \cite{MM}. 

\begin{definition}
\label{def:alexprinc}
Let $G$ be a finitely generated group. 
We say that the Alexander ideal $\E_1(A_G)$ is 
{\em almost principal} if there exists an integer 
$d\ge 0$ such that 
\[
I^d\cdot (\Delta^G)\subset \E_1(A_G),
\]
over $\C$, where $I$ denotes the augmentation ideal of the 
group ring $\Z[G_{\abf}]$.
\end{definition}

For this class of groups, the reduced Alexander polynomial
determines $\V_1(G)\cap \T_G^0$, since 
$\V_1(G)\cap \T_G^0= V(\Delta^G)$, away from $1$, 
as noted in Proposition \ref{prop:echar}.  There is an abundance 
of interesting examples of this kind.
 
Recall that the {\em deficiency} of a finitely presented group 
$G$, written $\df(G)$, is the supremum of the difference between 
the number of generators and the number of relators, taken 
over all finite presentations of $G$. 

\begin{lemma}
\label{lem:b1=1 or pos def}
If $b_1(G)=1$, or $\df (G)>0$, then $\E_1(A_G)$ is 
almost principal.
\end{lemma}

\begin{proof}
If $b_1(G)=1$, then $\Lambda_1(\C)$ is principal, and so 
$\E_1(A_G)= (\Delta^G)$. If $b_1(G)\ge 2$ and $\df (G)>0$, 
then  $\E_1(A_G)= I\cdot (\Delta^G)$, by Theorems 6.1 and 
6.3 in \cite{EN}.
\end{proof}

Another class of groups with this property occurs in 
low-dimensional topology. 

\begin{lemma}
\label{lem:3mfd}
Let $M$ be a compact, connected $3$-manifold, and let 
$G=\pi_1(M)$. If either $\partial M\ne \emptyset$ and 
$\chi (\partial M)=0$, or $\partial M=\emptyset$ and $M$ 
is orientable, then $\E_1(A_G)$ is almost principal.
\end{lemma}

\begin{proof} 
In the first case, $\df(G)> 0$, by \cite[Lemma 6.2]{EN}, and so 
the result follows from Lemma \ref{lem:b1=1 or pos def}.  In the 
second case, $I^2\cdot (\Delta^G)\subset \E_1(A_G)$, as shown 
in \cite[Theorem 5.1]{MM}.
\end{proof}

\subsection{Multiplicity bounds for $b_1(G, \C_{\rho})$}
\label{subsec:aprinc bounds}

We are now ready to state our result concerning multiplicity bounds 
for the twisted Betti numbers for the class of groups delineated 
above. 

For a character $\rho\in (\C^*)^n$, and a Laurent polynomial 
$f\in \Lambda_n$, denote by $\nu_{\rho}(f)$ the order of vanishing 
of the germ of $f$ at $\rho$.  

\begin{theorem}
\label{thm:uppergral}
Let $G$ be a finitely generated group, and let  
$\Delta^G \doteq_{\C} f_1^{\mu_1}\cdots f_s^{\mu_s}$ be 
the decomposition into irreducible factors of its 
Alexander polynomial.  Assume the Alexander ideal of 
$G$ is almost principal. 
If $\rho\in \T_G^0\setminus \{1\}$, then
\begin{equation}
\label{eq:bd}
\dim_{\C} H_1(G; \C_{\rho})\le \sum_{j=1}^s \mu_j\cdot \nu_{\rho}(f_j).
\end{equation}
\end{theorem}

\begin{proof}
As before, we may assume $G$ admits a finite presentation, say 
$G= \langle x_1, \dots, x_m \mid r_1, \dots, r_h \rangle$. 
Set $n=b_1(G)$, and consider the change of rings 
$\psi\colon \Lambda_n \to \OO_n$ given by restriction 
of Laurent polynomials near a fixed $\rho\in \T_G^0 \setminus \{ 1\}$. 
Let $\AA_G^{\psi} \colon \OO_n^h \to \OO_n^m$ 
be the corresponding `local' Alexander matrix, and 
$\AA_G(\rho)\colon \C^h \to \C^m$ its evaluation at $\rho$. 
Set $r=\rank \AA_G(\rho)$ and $b=b_1(G, \rho)$. 
We then have
\begin{equation}
\label{eq:prel}
\dim_{\C} \coker \AA_G(\rho)= 
m-r= b+1 ,
\end{equation}
since $\rho\ne 1$. 

Setting $Z=\ker \AA_G(\rho)$,  
we may choose vector space decompositions, $\C^h= Z\oplus \C^r$ 
and $\C^m= N\oplus \C^r$, such that $\AA_G(\rho)$ takes 
the diagonal form $0\oplus \id$. In the corresponding 
decomposition over $\OO_n$, the $r\times r$ submatrix 
of $\AA_G^{\psi}$ is invertible.  Using \eqref{eq:prel}, we see that 
\begin{equation}
\label{eq:min}
\coker \AA_G^{\psi}= 
\coker \big( D\colon \OO_n^{h-r}\rightarrow \OO_n^{b+1} \big) ,
\end{equation}
where $D$ is equivalent to the zero matrix, modulo  
the maximal ideal $\m$ of $\OO_n$.

By assumption, the Alexander ideal of $G$ is almost principal, 
i.e., $I^d\cdot (\Delta^G)\subset \E_1(A_G)$, for some $d$. 
Since $\rho\ne 1$, we obtain by base change
\[
\psi (\Delta^G)\in \E_1(\coker \AA_G^{\psi}). 
\]
Using \eqref{eq:min}, we deduce that $\psi (\Delta^G)\in \m^b$. 
Hence, $b\le \nu_{\rho}(\Delta^G)=\sum_j \mu_j \cdot \nu_{\rho}(f_j)$, 
as asserted.
\end{proof}

\begin{corollary}
\label{cor:deltadetv}
Assume in Theorem \ref{thm:uppergral} that $G_{\ab}$ is 
torsion-free. If the upper bound \eqref{eq:bd} is attained for 
every  $\rho\in \T_G^0\setminus \{1\}$, then the Alexander 
polynomial $\Delta^G$ determines the characteristic 
varieties $\V_k(G)$, for all $k\ge 1$. 
\end{corollary}

\subsection{Discussion} 
\label{susbec:discuss uppergral}

As the next example shows, the hypothesis that the character 
$\rho$ be non-trivial is essential for Theorem \ref{thm:uppergral} 
to hold. 

\begin{example}
\label{rem:non-trivial system}
Let $G=\Z\times F_{n-1}$, where $n\ge 3$.   Note that $\df(G)=1$, 
and so, by Lemma \ref{lem:b1=1 or pos def}, the Alexander ideal 
of $G$ is almost principal. Clearly, $b_1(G, 1)=n$, 
while the upper bound from Theorem \ref{thm:uppergral} 
equals $(n-2)\cdot \nu_1(t_1\cdots t_n -1)=n-2$. 
\end{example}

The hypothesis on the Alexander ideal being almost principal 
is also needed.

\begin{example}
\label{ex:notprinc}
The method discussed in Example \ref{ex:alexvsgen} 
produces a commutator-relators group, 
$G=\langle x_1, x_2, x_3\mid r_1, r_2, r_3\rangle$, 
with Alexander matrix
\[
\AA_G\begin{pmatrix}
(x_2-1) \Phi & (1-x_1) \Phi & 0\\
0 & (x_3-1) \Phi & (1-x_2) \Phi\\
(x_3-1) \Psi & 0 & (1-x_1) \Psi
\end{pmatrix},
\]
where $\Phi =x_1x_3 -1$ and $\Psi = x_1x_2 -1$. 
As is readily seen, $\Delta^G \doteq \Phi$. For the nontrivial local system 
$\rho= (1,-1,1)$, we have $b_1(G, \rho)=2$; on the other hand, the 
corresponding upper bound from \eqref{eq:bd} is $1$. 
\end{example}

Inequality \eqref{eq:bd}  may well be strict, as the next example shows. 

\begin{example}
\label{ex:strict ineq}
The groups $G_k$, with $k>1$, from Example \ref{ex:gensharp} 
meet the requirements from Theorem \ref{thm:uppergral}. If 
$\rho\in V(f)\setminus \{ 1\}$, then clearly $b_1(G_k, \rho)=1< k \nu_{\rho}(f)$. 
\end{example}

Finally, let us point out that inequality \eqref{eq:bd} is sharp.

\begin{example}
\label{ex:gralsharp}
Let $G=\Z\times F_{n-1}$ be the group from Example \ref{ex:pencil}, 
with Alexander polynomial $\Delta^G=f^{n-2}$, where 
$f=t_1\cdots t_n -1$.  Clearly, $G$ satisfies the assumptions 
of Theorem \ref{thm:uppergral}.  If $\rho\in V(f)\setminus \{ 1\}$, 
then $b_1(G, \rho)= (n-2)\cdot \nu_{\rho}(f)=n-2$, i.e., the upper
bound \eqref{eq:bd} is attained.
\end{example}

\subsection{The bounds in a special case}
\label{subsec:single multbound}

When $\Delta^G \doteq_{\C} \prod_j f_j^{\mu_j}$ has a single essential 
variable, Theorem \ref{thm:uppergral} takes the following simpler form.

\begin{corollary}
\label{cor:upperqp}
Suppose the Alexander ideal of $G$ is almost principal, 
$\Delta^G$ has a single essential variable, 
and $\Delta^G \not\doteq \const$. Then 
\begin{enumerate}
\item \label{c1}
The intersection $\V_1(G)\cap \T_G^0 \setminus \{ 1 \}$ equals 
the disjoint union $\coprod_j V(f_j)\setminus \{ 1\}$. 

\item \label{c2}
If $1\ne \rho\in V(f_j)$, then $1\le b_1(G, \rho)\le \mu_j$.
\end{enumerate}
\end{corollary}

\begin{proof}
Part \eqref{c1} follows by putting together Definition \ref{def:essvar} 
and the remark made after Definition \ref{def:alexprinc}. As for 
Part \eqref{c2}, use the obvious fact that $\nu_{\rho}(t_1-z)=1$, 
for any $\rho\in V(t_1-z)$. 
\end{proof}

\begin{remark}
\label{rk:trivmono}
Let $G$ be a group as in Corollary \ref{cor:upperqp}. 
Pick any nontrivial element $\rho\in V(f_j)$, 
and let $N$ be a local transversal slice to $V(f_j)$ at $\rho$. 
In this situation, the arguments from the proof of 
Theorem \ref{thm:uppergen} work not only for generic 
elements of $V(f_j)$, but also for arbitrary $\rho$ in  
$\T_G^0 \setminus \{ 1 \}$, giving the following information.

In \eqref{eq:cokdec}, $r=1$ and $b_1 (G,0)=\ell$. Due to the fact 
that $\E_1(A_G)$ is almost principal, and $\rho \ne 1$, we have 
$(\tau(\Delta^G))= \E_1(\coker \AA_G^{\tau})= 
(z^{\nu_1+\cdots +\nu_{\ell}})$.  This implies
\[
\mu_j= \nu_0\Delta^G(z, 0)= 
\sum_{k=1}^{\ell} \nu_k \ge \ell= b_1(G, \rho) .
\]
Therefore, the equality $b_1(G, \rho)=\mu_j$ is equivalent to
\begin{equation}
\label{eq:1jordan}
\nu_1=\cdots =\nu_{\ell}=1 ,
\end{equation}
or, in invariant form, $z\cdot \big( \OO_1(\rho, N)\otimes_{\C[G_{\abf}]} 
B_G \otimes \C\big)= 0$, 
where $\OO_1(\rho, N)$ denotes $\OO_1$, with module structure 
given by restriction of Laurent polynomials to $N$.
\end{remark}

We shall see in Theorem \ref{thm:seifert cv} that conditions \eqref{eq:1jordan} 
are satisfied by  Seifert links.  
On the other hand, these conditions are not satisfied in general, 
as the next example shows (see also Examples \ref{ex:rolfsen} 
and \ref{ex:notss} below). 

\begin{example}
\label{ex:one more}
Let $G=\langle x_1, x_2\mid r_1, r_2\rangle$ be a 
commutator-relators group, with $b_1(G)=2$ and Alexander matrix
\[
\AA_G =
\begin{pmatrix}
(x_1-1)(x_2-1)\alpha & -(x_1-1)^2 \alpha\\
(x_2-1)^2 \alpha & -(x_1-1)(x_2-1)\alpha
\end{pmatrix},
\]
where $\alpha= (x_1x_2 +1)^k$, and $k>1$. In this case, 
$\Delta^G(x_1,x_2) \doteq \alpha$; in particular, $\Delta^G$ 
has a single essential variable.  Moreover, 
$I^2\cdot (\Delta^G)\subset \E_1(A_G)$, and so the 
Alexander ideal is almost principal. Now, for any 
$\rho\in V(x_1x_2 +1)\subset \T_G \setminus \{ 1\}$, we 
see that $b_1(G, \rho)= 1$, even though the multiplicity 
of the factor $x_1x_2 +1$ of $\Delta^G$ is $k>1$. 
\end{example}

\subsection{Monodromy}
\label{subsec:mono}

For the rest of this section, we assume $b_1(G)=1$. 
In this case, the torsion-free abelianization map is 
$\phi_{\abf}\colon G\to \Z$. Identify the group ring $\Z\Z$ 
with $\Lambda_1= \Z[t^{\pm 1}]$, 
and note that $I_{\Z}=(t-1)\Lambda_1$ is a free $\Lambda_1$-module 
of rank $1$. The Crowell exact sequence \eqref{eq:crow2} then 
yields a decomposition $A_G\cong B_G \oplus I_{\Z}$, from which 
we infer
\begin{equation}
\label{eq:e1ab}
\E_1(A_G)= \E_0(B_G).
\end{equation}

Obviously, the Alexander polynomial $\Delta^G$ has 
a single variable.  Since $\Lambda:=\Lambda_1\otimes \C\C[t^{\pm 1}]$ 
is a PID, the ideal $\E_1(A_G)\otimes \C$ is 
principal (generated by $\Delta^G$).  Hence, 
Corollary \ref{cor:upperqp} applies, provided 
$\Delta^G \not\doteq \const$. In view of \eqref{eq:e1ab}, this 
condition is equivalent to $B_G\otimes \C$ being a non-zero, 
torsion $\Lambda$-module, in which case 
\begin{equation}
\label{eq:princdec}
B_G\otimes \C = \bigoplus_{j=1}^{s}\bigoplus_{k\ge 1} 
\big( \Lambda/ (t-z_{j})^{k} \Lambda \big)^{e_k(z_{j})} ,
\end{equation}
where the sum over $k$ is finite, $z_{1},\dots ,z_{s}$ are distinct 
elements in $\C\setminus \{ 0,1 \}$, and $e_k(z_{j})\ge 1$, for all 
$k$ and $j$.

It follows from \eqref{eq:e1ab} and \eqref{eq:princdec} that the 
Alexander polynomial factors as 
\begin{equation}
\label{eq:delta z}
\Delta^G \doteq_{\C} (t- z_{1})^{\mu_1}\cdots (t- z_{s})^{\mu_{s}},
\end{equation}
with
\begin{equation}
\label{eq:princalex}
\mu_{j} = \sum_{k\ge 1} k\, e_k(z_{j}). 
\end{equation}

Every complex number $z\ne 0$ defines a character 
$z\in \T_{G}^0$, and a local system $\C_{z}$.  The next result 
relates the twisted Betti ranks corresponding to the roots 
$z_j$ of the Alexander polynomial to the exponents 
$e_k(z_j)$ appearing in the decomposition \eqref{eq:princdec} 
of the complexified Alexander invariant. 

\begin{prop}
\label{prop:princtest}
Let $G$ be a group with $b_1(G)=1$ and 
$\Delta^G \not\doteq \const$.  Then, for $j=1,\dots, s$, 
\[
b_1(G, z_j)=\mu_j  \Longleftrightarrow 
e_k(z_j)=0,\ \forall k>1.
\]
\end{prop}

\begin{proof}
Since $\Lambda$ is a PID and $z_j\ne 1$, we have 
\[
b_1(G, z_j)=\dim_{\C} (\C_{z_j}\otimes_{\Lambda} B_G\otimes \C)\sum_{k\ge 1} e_k(z_j).
\]
Comparing this equality with \eqref{eq:princalex} establishes the claim.
\end{proof}

The conditions $e_k(z_j)=0$ for $k>1$ represent a delicate 
restriction on the size of the Jordan blocks of a presentation 
matrix for the Alexander invariant. These conditions are not 
always satisfied, even for  knots in $S^3$.  

\begin{example}
\label{ex:rolfsen}
By a classical result of Seifert (see \cite[Theorem 7C.5]{Ro}), 
given any Laurent polynomial $f\in \Lambda_1=\Z[t^{\pm 1}]$ with 
$f(1)=\pm 1$ and $f(t)\doteq f(t^{-1})$, there exists  
a knot in $S^3$ for which the Alexander invariant of 
the complement is isomorphic to $\Lambda_1/(f)$. 

So fix such a polynomial, say $f(t)=t^2-t+1$, 
and for each $k>1$, construct a knot $K\subset S^3$ with 
group $G=\pi_1(S^3\setminus K)$ and Alexander 
invariant $B_G=\Lambda_1/(f^k)$.  Clearly, all roots $z_j$ of
$\Delta^G= f^k$ are multiple roots, and so $b_1(G,z_j)<\mu_j$.
\end{example}

\subsection{Fibrations over the circle}
\label{subsec:fibrations}
Let $M$ be a compact, connected manifold (possibly with 
boundary), which admits a locally trivial fibration over the circle,
$F \hookrightarrow M \xrightarrow{p} S^1$. Let $G=\pi_1(M)$, 
and denote by $h\colon F\to F$ the monodromy of the fibration.  

Suppose the matrix of $h_*\colon H_1(F;\C)\to H_1(F;\C)$ does 
not have $1$ as an eigenvalue.  By the Wang sequence of the 
fibration, the map $p_* \colon H_1(M; \C)\to H_1(S^1; \C)$ 
is an isomorphism, and so $b_1(G)=1$. If $F$ is connected, 
the Alexander polynomial of $G$ equals (up to units in 
$\Lambda$) the characteristic polynomial of $h_*$; in 
particular, $\Delta^G \not\doteq \const$, when $b_1(F)>0$. 
Hence, Proposition \ref{prop:princtest} applies, and we obtain the following 
corollary. 

\begin{corollary}
\label{cor:monotest}
In the above setup, write 
$\Delta^G \doteq_{\C} \prod (t- z_j)^{\mu_j}$, 
as in \eqref{eq:delta z}. The following are equivalent:
\begin{enumerate}
\item $b_1(G, z_j)=\mu_j$, for all $j$. 
\item The algebraic monodromy, $h_*\colon H_1(F;\C)\to H_1(F;\C)$, 
is semisimple.
\end{enumerate}
\end{corollary}

\begin{proof}
Since $p_* \colon H_1(M; \Z)\to H_1(S^1; \Z)$ is onto, $h_*$ may be
read off from \eqref{eq:princdec}.
\end{proof}

\begin{example}
\label{ex:notss}
As is well-known, any symplectic matrix $A\in\Sp(2g, \Z)$ 
may be realized as $h_*\colon H_1(F; \Z)\to H_1(F; \Z)$, for 
some homeomorphism $h\colon F\to F$ of a compact, orientable 
Riemann surface $F$ of genus $g$, see \cite[Theorem~N13]{MKS}.  
The mapping torus of this homeomorphism, 
call it $M$, fibers over the circle, with monodromy $h$. 
If $1$ is not an eigenvalue of $A$, we are in the 
setup from Corollary \ref{cor:monotest}. If, moreover, 
$A$ is not diagonalizable, then the upper bounds on the 
twisted Betti numbers of $G=\pi_1(M)$ are not attained.

As a concrete example, take $A=\left( \begin{smallmatrix} -1 & 1\\
0 & -1 \end{smallmatrix} \right)$.  Then $M$ is a torus bundle 
over the circle, with fundamental group 
\[
G=  \langle x_1,x_2,x_3 \mid x_1 x_2 =x_2 x_1, \:
x_3^{-1} x_1 x_3 = x_1^{-1}, \: 
x_3^{-1} x_2 x_3 = x_1 x_2^{-1}
\rangle
\] 
and Alexander polynomial $\Delta^G=(1+t)^2$.  
For the local system defined by the eigenvalue $z_1=-1$, 
we have $b_1(G,z_1)=1$, which is less than the multiplicity 
$\mu_1=2$.
\end{example}

\section{Seifert links}
\label{sect:seifert}

In this section, we examine the class of Seifert links, considered by 
Eisenbud and Neumann in \cite{EN}. For such links, the Alexander 
polynomial of the link group $G$ determines {\em all} the characteristic 
varieties $\V_k(G)$, $k\ge 1$. 

\subsection{Links in homology $3$-spheres}
\label{subsect:links}
Let $\Sigma^3$ be a compact, smooth $3$-manifold with  
$H_*(\Sigma^3;\Z)\cong H_*(S^3;\Z)$.  A $q$-component 
link in $\Sigma^3$ is a collection of disjoint circles, 
$S_1, \dots, S_q$ ($q\ge 1$), smoothly embedded in $\Sigma^3$.  
We shall fix orientations on both $\Sigma^3$ and the link 
components $S_i$, and denote the resulting oriented link by 
\[
L= (\Sigma, S_1\cup \cdots \cup S_q). 
\]

The link exterior is the closed $3$-manifold 
$M_L=\Sigma^3 \setminus N_L$, where $N_L$ is an open 
tubular neighborhood of $L$. Clearly, $M_L$ has the same 
homotopy type as the link complement, $\Sigma^3 \setminus L$.  
The link group, $G_L:=\pi_1(M_L)$, is a finitely presented group, 
of positive deficiency.  Its abelianization, $(G_L)_{\ab} =\Z^q$, comes 
endowed with a canonical basis, given by oriented meridians around 
the link components, $e_1,\dots ,e_q$; in particular, 
$\T_{G_L}= \T_{G_L}^0$. By either Lemma \ref{lem:b1=1 or pos def} 
or Lemma~\ref{lem:3mfd}, the Alexander ideal of $G_L$ is 
almost principal. 

\subsection{The Eisenbud-Neumann calculus}
\label{ss:Seifert links}
If the link exterior $M_L$ admits a Seifert fibration, then $L$ is 
called a {\em Seifert link}.  Such links are completely classified. 
In this subsection, we consider only {\em positive} Seifert links, 
i.e., links of the 
form $L=(\Sigma(k_1,\dots, k_n), S_1\cup \cdots \cup S_q)$, 
with $k_j\ge 1$ pairwise coprime integers, $n\ge 3$ 
and $n\ge q\ge 1$, in the notation from \cite[Proposition 7.3]{EN}. 
Such links are conveniently represented by 
`splice diagrams'.  The diagram of the above link $L$ 
has one positive {\em node}, connected to {\em arrowhead vertices} 
$\{ v_1,\dots, v_q \}$ and {\em boundary vertices} 
$\{ v_{q+1},\dots, v_n \}$, with the vertices $\{ v_1,\dots, v_n \}$
being labeled by the integer weights $\{ k_1,\dots, k_n\}$; 
see \cite[p.~69]{EN}.

This restriction is motivated by the fact that fundamental groups 
of non-positive Seifert links are very simple: they are either 
free groups, or groups of the form $\Z\times F_{n-1}$, with 
$n\ge 2$;  see \cite[Proposition 7.3]{EN}. In all these cases,
the Alexander polynomial has a single essential variable.

Theorem 12.1 from \cite{EN} allows one to read off the Alexander 
polynomial $\Delta^L:=\Delta^{G_L}$ directly from the splice diagram. 
We assume for simplicity that $q\ge 2$. Then
\begin{equation} 
\label{eq2}
\Delta ^L(t_1,\dots,t_q)=\frac{[(t_1^{N_1}\cdots t_q^{N_q})^{N'}-1]^{n-2}}%
{\prod_{j=q+1}^{n} [(t_1^{N_1}\cdots t_q^{N_q})^{N'_j}-1]}, 
\end{equation}
where $N=k_1 \cdots k_q$, $N_j=N/k_j$ for $1\le j \le q$, 
$N'=k_{q+1} \cdots k_n$, and $N'_j=N'/k_j$ for $q+1\le j \le n$.
We can assume that for $j>q$, we have $k_j=1$ iff $j>q+s$, 
with $s$ a positive integer. Consequently, $N_j'=N'$ for $j>q+s$. 
After simplification, formula \eqref{eq2} becomes
\begin{equation} 
\label{eq3}
\Delta ^L(t_1,\dots,t_q)=\frac{[(t_1^{N_1}\cdots t_q^{N_q})^{N'}-1]^{q+s-2}}%
{\prod_{j=q+1}^{q+s} [(t_1^{N_1}\cdots t_q^{N_q})^{N'_j}-1 ]} .
\end{equation}

Since the integers $k_j$ are pairwise coprime, it follows that 
$\gcd(N_1,\dots,N_q)=1$, and hence the equation
\begin{equation} 
\label{eq4}
 t_1^{N_1}\cdots t_q^{N_q}=1     
\end{equation}
defines a subtorus $\T' \subset \T:=(\C^*)^q$. In particular, $\Delta^L$ 
has a single essential variable, namely $u= t_1^{N_1}\cdots t_q^{N_q}$, 
and so Corollary \ref{cor:upperqp} applies.

\subsection{Seifert links and quasi-projectivity}
\label{ss:qp seifert}

It follows from the preceding subsections that 
Corollary \ref{cor:upperqp} actually applies to all Seifert 
links with non-constant Alexander polynomial.
Now we turn to a discussion of some analytic aspects of Seifert links.

Let $(X,0)$ be a complex quasi-homogeneous normal surface 
singularity.   Then the surface $X^*=X \setminus \{0\}$ is smooth, 
and admits a $\C^*$-action with finite isotropy groups $\C^*_x$. 
These isotropy groups are trivial, except for 
those corresponding to finitely many orbits associated 
to some points $p_1, \dots, p_s$ in $X^*$; let 
$k_j$ be the order of $\C^*_{p_j}$. 

The quotient $X^*/\C^*$ is a smooth projective curve. 
For any finite subset $B\subset X^*/\C^*$, there is a 
surjective map $f\colon M \to S$, induced by the quotient map
$f_0\colon X^* \to X^*/\C^*$, where $S$ is the complement 
of $B$ in $X^*/\C^*$, and $M=f_0^{-1}(S)$.  Note that 
$f_*\colon  H_1(M) \to H_1(S)$ is surjective, since the 
fibers of $f$ are connected.  

In addition, the curve $X^*/\C^*$ is rational if and only if 
the link $L(X)$ of the singularity $(X,0)$ is a $\Q$-homology 
sphere (use Corollary (3.7) on p.~53 and Theorem (4.21) 
on p.~66 in \cite{D1}). In particular, if the link $L(X)$ of the 
singularity $(X,0)$ is a $\Z$-homology 
sphere, then $ H_1(M)=\Z^q$ where $q=\abs{B}$, and 
a basis is provided by small loops $\gamma_b$ around 
the fiber $F_b=f_0^{-1}(b)$ for $b \in B$. We have 
$f_*(\gamma_b)=k_b \delta_b$, with $k_b$ the order 
of the isotropy groups of points $x$ such that $f_0(x)=b$, 
and $\delta_b$ a small loop around $b \in \PP^1$.
The set of critical values of the map 
$f_0\colon X^* \to X^*/\C^*$ is exactly $B_0=\{y_1,\dots ,y_s\}$, 
with $y_j=f_0(p_j)$, and each fiber $F_j=f_0^{-1}(y_j)$ 
is smooth (isomorphic to $\C^*$), but of multiplicity $k_j>1$.
Writing down the map $f_{0*}$ and using its surjectivity, we 
deduce that the integers $k_1,\dots, k_s$ are pairwise coprime.

\begin{remark} 
\label{rk1}
Let $(X,0)$ be the germ of an isolated complex surface singularity, 
such that the corresponding link $L_X$ is an integral homology sphere. 
Let $(Y,0)$ be a curve singularity on $(X,0)$. Then using the conic structure
of analytic sets, we see that the local complement $X \setminus Y$, with 
$X$ and $Y$ Milnor representatives of the singularities $(X,0)$ and $(Y,0)$, 
respectively, has the same homotopy type as the link complement 
$M=L_X \setminus L_Y$, where $L_Y$ denotes the link of $Y$.

Moreover, if $(X,0)$ and $(Y,0)$ are quasi-homogeneous singularities 
at the origin of some affine space $\C^N$, with respect to the same 
weights, then the local complement can be globalized, i.e., replaced 
by the smooth quasi-projective variety $X \setminus Y$, where 
$X$ and $Y$ are this time affine varieties representing the germs
$(X,0)$ and $(Y,0)$, respectively.

A description along these lines of all Seifert links $L$
is given in \cite{EN}, p.~62. It follows that the link group $G_L$
is quasi-projective.
\end{remark}

\subsection{Characteristic varieties of Seifert links}
\label{ss:cv seifert}

We are now in position to describe all characteristic 
varieties of Seifert links, solely in terms of their Alexander 
polynomials.
 
\begin{theorem} 
\label{thm:seifert cv}
Let $G_L$ be the group of the Seifert link
$L=(\Sigma(k_1,\dots, k_n), S_1\cup \cdots \cup S_q)$, with 
$k_i\ge 0$, $n\ge 3$ and $n\ge q\ge 1$. Assume 
$\Delta^L \not\doteq \const$, that is, $G_L$ is not a free group 
or $\Z^2$.  Let $D=m_1D_1+\cdots+m_pD_p$ ($m_j>0$) be the 
effective divisor defined by $\Delta ^L=0$ in 
$\T_{G_L}$. Then, for all $j=1,\dots,p$, and all non-trivial 
characters $\rho \in D_j$, 
\[
\dim_{\C} H_1(G_L;\C_{\rho})= m_j  .
\]
\end{theorem}

\begin{proof}
Since $H_1(G;\C_{\rho})$ and $H^1(G;\C_{\rho})$ are dual 
vector spaces, with $\C_{\rho}$ viewed as a right $\C{G}$-module 
for homology and as a left $\C{G}$-module for cohomology, 
we may freely switch from homology to cohomology and back.

In the case of non-positive Seifert links, we only need to examine 
the groups from Example \ref{ex:pencil}, with $n\ge 3$. This was 
done in Example \ref{ex:gralsharp}.  Exteriors of positive Seifert 
knots fiber over $S^1$, with connected fibers and with 
finite-order geometric monodromy $h$; see \cite[Lemma 11.4]{EN}. 
In this case, the result follows from Corollary \ref{cor:monotest}. 
Thus, we may suppose from now on that $k_i\ge 1$, for all $i$, 
and $q\ge 2$.

Using the analytic description of the Seifert link 
$L=(\Sigma (k_1,\dots,k_n), S_1 \cup\cdots\cup S_q)$ 
recalled in Remark \ref{rk1} 
and the notation introduced in Section \ref{ss:qp seifert}, 
we see that the link complement,  
$M(L)=\Sigma (k_1,\dots,k_n) \setminus \{ S_1 \cup\cdots\cup S_q \}$, 
has the homotopy type of the surface $M$ obtained from the surface 
singularity $X$ by deleting the orbits (regular for $k_j=1$ and singular 
for $k_j>1$) corresponding to the knots $S_1,\dots, S_q$. 
In other words, we have a finite set $B \subset \PP^1$ with 
$\abs{B}=q$ and an admissible mapping $f\colon M \to 
S=\PP^1 \setminus B$, exactly as in \S \ref{ss:qp seifert}.

Using the basis $\{\gamma_b\}_{b\in B}$ of $H_1(M)$---%
respectively, the system of generators $\{\delta_b\}_{b\in B}$---%
described in \S \ref{ss:qp seifert}, it follows that $f$ induces an 
embedding of tori, $f^*\colon\T(S) \to \T(M)$, whose image is 
exactly the subtorus $\T'$ defined in \eqref{eq4}. 
Thus, we have an isomorphism
\[
\T(M)/f^*(\T(S)) \isom \C^*, 
\]
given by $(t_1,\dots,t_q) \mapsto  t_1^{N_1}\cdots t_q^{N_q}$. 
On the other hand, let $i\colon F \to M$ be the inclusion of a 
generic (i.e., non-multiple) fiber of $f$.  We know that $F=\C^*$, 
hence we get a restriction morphism
\[
i^*\colon\T(M) \to \T(F)=\C^*.
\]

One can calculate the image of $i_*\colon H_1(F) \to H_1(M)$, 
as follows. The set of critical values $C(f)$ of $f\colon M \to S$
consists precisely of the $s$ points in $S$ corresponding to the 
knots $S_j$ (which are not used in the link, i.e., for $j>q$) with 
multiplicity $k_j>1$. Then, if we delete these critical values and 
the corresponding multiple fibers, we obtain a locally trivial fibration 
$f'\colon M' \to S'$. If $i'\colon F \to M'$ denotes the inclusion,
it is easy to determine $\im(i'_*)=\ker(f'_*)$ and then one obtains  
$\im(i_*)$ using the natural projection $H_1(M') \to H_1(M)$ 
(see \cite{D06} for a similar computation). If $\sigma \in H_1(F)$ 
denotes the natural generator, it follows that
$i_*(\sigma)=N'(N_1\gamma_1+\cdots +N_q\gamma_q)$. 
Hence,
\begin{equation} 
\label{eq4.5}     
i^*(t_1,\dots ,t_q)=(t_1^{N_1}\cdots t_q^{N_q})^{N'}.
\end{equation}

Using \cite{D06}, it follows that the  irreducible components of 
$\V_1(M)=\V_1(G_L)$ associated to the admissible map 
$f\colon M \to S$ are parametrized by 
\[
\widehat T(f)=\frac{\ker i^*}{f^*(\T(S))}=\Gamma _{N'} ,
\]
for $q>2$, respectively by $\Gamma _{N'}\setminus \{ 1\}$, 
for $q=2$, where  $\Gamma _{N'}$ is the cyclic group of 
roots of unity of order $N'$.  By virtue of \eqref{eq3} and 
Corollary \ref{cor:deltacd1}\eqref{dc2}, these components 
coincide with the irreducible components of $V(\Delta^L)$.

Now, for a character $\rho \in \widehat T(f)$, it follows 
from \cite{D06} that
\begin{equation} 
\label{eq5}     
\dim H^1(M;\C_{\rho'})\ge -\chi(S)+\abs{\Sigma (R^0f_*(\C_{\rho}))} ,
\end{equation}
for any $\rho'$ in the corresponding irreducible component 
$W_{f,\rho}$, provided the right hand side is strictly positive. 
Moreover, equality holds with finitely many exceptions.

Using formula \eqref{eq3} and Corollary \ref{cor:upperqp}, 
we obtain
\begin{equation} 
\label{eq6}     
\dim H^1(M;\C_{\rho'})\le q+s-2-\abs{I(\alpha(\rho))} ,
\end{equation}
for $\rho'\ne 1$,
where $\alpha(\rho)= \rho_1^{N_1}\cdots \rho_q^{N_q}$ and 
$I(\alpha)=\{j \mid q<j\le q+s \text{ and } \alpha^{N_j'}=1\}$.

On the other hand,  $\chi(S)=2-q$. If $c \in C(f)$ 
and $T(F_c)=f^{-1}(D_c)$ is a small tubular neighborhood 
of the multiple fiber $F_c$, with $D_c$ a small open disk  
centered at $c$, then the inclusion $F_c \to T(F_c)$ is a 
homotopy equivalence; hence, $H_1(F_c)=H_1(T(F_c))$.
The inclusion $F \to T(F_c)$ induces then a morphism $
H_1(F) \to H_1(T(F_c))$ which is just multiplication by 
$k_c\colon \Z \to \Z$. Combining this fact with 
formula \eqref{eq4.5}, it follows that the inclusion 
$i_c\colon F_c\to M$ induces a morphism 
$i_c^*\colon \T(M) \to \T(F_c)=\C^*$ given by
\begin{equation} 
\label{eq6.5}     
i_c^*(t_1,\dots ,t_q)=(t_1^{N_1}\cdots t_q^{N_q})^{N'_c}.
\end{equation}

Now recall from \cite{D06} that $c \in \Sigma (R^0f_*(\C_{\rho}))$ 
if and only if $H^0(T(F_c);\C_{\rho})=0$. This happens exactly 
when the restriction of the local system $\C_{\rho}$ to $T(F_c)$---or, 
equivalently, to $F_c$---is non-trivial. By formula \eqref{eq6.5}, 
this is the same as $c \notin I(\alpha(\rho))$, i.e.,
\begin{equation} 
\label{eq:sigmar0} 
\abs{\Sigma (R^0f_*(\C_{\rho}))}=s-\abs{I(\alpha(\rho))}.
\end{equation}
This establishes the reverse inequality in \eqref{eq6}. Since the right hand side
from \eqref{eq6} is precisely the multiplicity from the assertion to be proved,
we are done.
\end{proof}

\begin{example} 
\label{ex4}
Let $p$ and $q$ be two positive integers such that $(p,q)=1$ and 
let $n \ge 2$.  Let $L$ be the algebraic link in $S^3$ associated 
to the plane curve singularity
\[
(C,0)\colon x^{pn}-y^{qn}=0\, .
\]
In terms of Seifert's notation, this is the link
\[
L=(\Sigma (1,\dots ,1,p,q), S_1 \cup\cdots \cup S_n) ,
\]
with $1$ occurring $n$ times. Formula \eqref{eq2} becomes
\[
\Delta ^L(t_1,\dots,t_n)=\frac{[(t_1\cdots t_n)^{pq}-1]^n}%
{[(t_1\cdots t_n)^{p}-1][(t_1\cdots t_n)^{q}-1]} .
\]

The irreducible components of the divisor $D$ in $\T ^n$ 
defined by $\Delta ^L(t_1,\dots,t_n)=0$ are parametrized by 
the elements of the group $\Gamma_p \times \Gamma_q$, 
for $n\ge 3$, respectively by the nontrivial elements, for $n=2$.
(Here $\Gamma _k$ denotes the multiplicative group of $k$-th roots 
of unity.)  More precisely, for $(a,b)\in \Gamma_p \times \Gamma_q$, 
the corresponding irreducible component
\begin{equation}
\label{eqd}
D_{a,b}\colon t_1\cdots t_n -a\cdot b=0
\end{equation}
(a subtorus through $1$ if $a\cdot b=1$, and a translated subtorus if 
$a\cdot b \ne 1$) occurs with multiplicity 
$m(a,b)=n-2+\abs{\{a-1,b-1\}\cap \C^*}$.

It follows that the characteristic variety $\V_1(G_L)$ has irreducible 
components given by \eqref{eqd}, with $(a,b)\in \Gamma_p \times
\Gamma_q$, for $n\ge 3$; if $n=2$, one needs to replace 
$D_{1,1}$ by $\{ 1\}$.
\end{example}

\begin{ack}
This work was started while the authors visited the Abdus Salam 
International Centre for Theoretical Physics in Fall, 2006. We thank 
ICTP for its support and excellent facilities, and L\^e D\~ung Tr\'ang 
for a helpful hint. 
\end{ack}

\newcommand{\arxiv}[1]
{\texttt{\href{http://arxiv.org/abs/#1}{arxiv:#1}}}

\renewcommand{\MR}[1]
{\href{http://www.ams.org/mathscinet-getitem?mr=#1}{MR#1}}

\end{document}